\documentclass[3p]{elsarticle}

\title{$p$-Adic Lifting Problems and Derived Equivalences}
\date{}

\usepackage{hyperref}
\usepackage{macros}
\usepackage[english]{babel}
\usepackage{longtable}
\usepackage{bibgerm}
\usepackage{url}
\usepackage[all]{xy}
\usepackage{ctable}

\usepackage{mathrsfs}

\usepackage{setspace}
\usepackage{enumitem}
\usepackage[active]{srcltx}
\newtheorem{defi}{Definition}[section]
\newtheorem{remark}[defi]{Remark}

\newtheorem{thm}[defi]{Theorem}
\newtheorem{lemma}[defi]{Lemma}
\newtheorem{corollary}[defi]{Corollary}
\newtheorem{notation}[defi]{Notation}
\newtheorem{prop}[defi]{Proposition}

\begin{document}

\begin{frontmatter}
 \author{Florian Eisele}
\ead{florian.eisele@rwth-aachen.de}
\address{Lehrstuhl D f\"ur Mathematik, RWTH Aachen, Templergraben 64, 52062 Aachen, Germany}
\begin{abstract}
    For two derived equivalent $k$-algebras $\overline\Lambda$ and $\overline\Gamma$, we introduce 
    a correspondence between $\OO$-orders reducing to $\overline\Lambda$ and $\OO$-orders reducing to $\overline\Gamma$.
    We outline how this may be used to transfer properties like uniqueness (or non-existence) of a lift between 
    $\overline\Lambda$ and $\overline\Gamma$. As an application, we look at tame algebras of dihedral type with
    two simple modules, where, most notably, we are able to show that among those algebras only the algebras $\mathcal D^{\kappa,0}(2A)$ and
  $\mathcal D^{\kappa,0}(2B)$ can actually occur as basic algebras of blocks of group rings of finite groups.
\end{abstract}
\end{frontmatter}

\section{Introduction}

In this article we consider, roughly speaking, the following problem: Given a finite-dimensional $k$-algebra $\overline{\Lambda}$, 
how many $\OO$-orders $\Lambda$ are there with $k\otimes \Lambda \iso \overline\Lambda$? Here $k$ is a field of
characteristic $p>0$ and $\OO$ is a complete discrete valuation ring with residue field $k$. Moreover, the field of fractions $K$ of $\OO$ is assumed to be of characteristic zero. Of course, there are usually infinitely many such lifts $\Lambda$ of $\overline\Lambda$,
and we may want to impose further restrictions on $\Lambda$ to get a meaningful answer. 
Our main focus lies on the case where $\overline\Lambda$ is a block of $kG$ for some finite group $G$, and the restrictions imposed on $\Lambda$ should therefore be
known properties of the corresponding block of $\OO G$. The semisimplicity of $K\otimes \Lambda$ and the symmetry of $\Lambda$ are certainly the simplest of those
properties, but further properties may come from knowledge of the decomposition matrix and character values of $G$.
The ideal outcome, for any given block of a group algebra $kG$, would be that any $\OO$-order $\Lambda$ subject to
a certain set of conditions is isomorphic to the corresponding block of $\OO G$.

For some algebras we can tackle these sorts of questions directly, essentially using linear algebra. 
The method we devise in this paper is the transfer, at least to some extent, of answers to the above questions via
derived equivalences of $k$-algebras. The key point is that although, for two given derived equivalent  finite-dimensional $k$-algebras $\overline\Lambda$ and $\overline\Gamma$, the respective problems of lifting $\overline\Lambda$ to an $\OO$-order and lifting
$\overline\Gamma$ to an $\OO$-order  may appear to be of a different degree of difficulty from the point of view of elementary linear algebra, we will show that they are in fact essentially equivalent.
For this we associate to a two-sided tilting complex $X \in \mathcal D^b(\overline\Lambda^{\op}\otimes_k\overline\Gamma)$
a bijection
\begin{equation}
	\Phi_X: \Lifts(\overline\Lambda) \longrightarrow \Lifts(\overline\Gamma)
\end{equation}
where $\Lifts(\overline\Lambda)$ denotes the set of equivalence classes of pairs $(\Lambda,\varphi)$, $\Lambda$ being an $\OO$-order
and $\varphi: k\otimes\Lambda \stackrel{\sim}{\rightarrow}\overline\Lambda$ being an isomorphism (see Definition \ref{defi_lifts} for a precise definition). $\Lifts(\overline\Gamma)$ is defined in the same way. Of course some more work is needed to make this map useful, as, for instance,
the sets $\Lifts(\overline\Lambda)$ and $\Lifts(\overline\Gamma)$ usually contain many elements representing one and the same order.

The idea behind this map $\Phi_X$ can be  explained fairly easily in the case of a Morita-equivalence:
If $\overline\Lambda$ and $\overline\Gamma$ are two Morita-equivalent $k$-algebras, then there is some
invertible $\overline\Lambda$-$\overline\Gamma$-bimodule $X$ with inverse $X^{-1}$. We can restrict $X^{-1}$ to a projective right $\overline\Lambda$-module $\overline P$. The endomorphism ring of $\overline P$ can be identified with  $\overline\Gamma$. 
Now for any pair $(\Lambda,\varphi)\in\Lifts(\overline\Lambda)$ we can use $\varphi$ to turn $\overline P$ into a projective $k\otimes\Lambda$-module.
This projective $k\otimes\Lambda$-module will lift uniquely to a projective $\Lambda$-module $P$. The endomorphism
ring (let us call it $\Gamma$) of $P$ will then be a lift of $\overline\Gamma$. What is still missing at this point is an isomorphism
$k\otimes \Gamma \stackrel{\sim}{\rightarrow} \overline\Gamma$. Since $P$ may be construed as a  $\Gamma$-$\Lambda$-bimodule,
we just choose $\psi: k\otimes \Gamma \stackrel{\sim}{\rightarrow}\overline\Gamma$ so that when we turn $X^{-1}$ into a $k\otimes\Gamma$-$k\otimes\Lambda$-bimodule using $\varphi$ and $\psi$, it becomes isomorphic to $k\otimes P$. Note that this does in fact determine $\psi$ uniquely. $(\Gamma,\psi)$ will then be the lift of $\overline\Gamma$ we associated to $(\Lambda,\varphi)$.

A field of application are tame blocks of group algebras over $k$. Here the appendix of \cite{TameClass} gives a list of $k$-algebras
which may occur as basic algebras. The first question is what the corresponding blocks of group algebras over 
$\OO$ look like, and under which conditions two tame blocks of group algebras over $\OO$ are Morita-equivalent given that the
corresponding blocks of the group algebras over $k$ are Morita-equivalent. 
Concretely, we look at blocks with dihedral defect group and two simple modules. The upshot here is that two such blocks
are Morita-equivalent over $\OO$ if and only if their corresponding blocks defined over $k$ are Morita-equivalent and their centers are equal.
Lifts for these blocks are then determined explicitly in Theorem \ref{thm_lifts_dihedral_two_simple}.

We also narrow down which algebras given in Erdmann's list
may actually occur as blocks of group rings. The key point here is that blocks of group rings defined over $k$ possess a lift, namely the
corresponding block defined over $\OO$, which has all the well-known properties that blocks of group rings share. 
We show that such a lift does not exist for algebras of dihedral type with two simple modules and parameter
$c=1$, implying that only those algebras with parameter $c=0$ occur as basic algebras of blocks of group rings (Corollary \ref{corr_ceq0_impossible}). The last assertion was recently proved for principal blocks in \cite{BleherUnivDef}, but the general case was still open.

\section{Foundations and Notation}

In this section we recall some definitions and theorems, and state some corollaries for later use.
Throughout this article, $(K,\OO,k)$ will denote a $p$-modular system for some $p>0$, and we assume that $\OO$ is complete. 
We let $\pi$ be a uniformizer of $\OO$.

\begin{notation}
	If $A$ is a ring, we denote by $\modC_A$ the category of finitely generated right $A$-modules, and by
	$\projC_A$ the category of finitely generated projective right $A$-modules. By a ``module'' we will always
	mean a right module (unless we explicitly say otherwise). By $\mathcal C^b(\projC_A)$ we denote the category
	of bounded complexes over $\projC_A$, and by $\mathcal K^b(\projC_A)$ we denote the corresponding homotopy category.
	By $\mathcal D^b(A)=\mathcal D^b(\modC_A)$ (respectively $\mathcal D^-(A)=\mathcal D^-(\modC_A)$)
	we denote the bounded derived category of $A$ (respectively the right-bounded derived category of $A$).
	By ``$-\otimes^{\mathbb L}_A=$'' we will always denote a left-derived tensor product.
\end{notation}


\begin{notation}
	Let $R\in\{\OO,k\}$, and let  $A$, $B$ and $B'$ be $R$-algebras. Let $\alpha: B \rightarrow A, \beta: B' \rightarrow A$ be $R$-algebra homomorphisms. Then define $_\alpha A _\beta$ to be the $B$-$B'$-bimodule
	which is (as a set) equal to $A$ with the action
	\begin{equation}
		B\times A \times B' \longrightarrow A:\ (b, x,b') \mapsto \alpha(b)\cdot x \cdot \beta(b') 
	\end{equation}
\end{notation}

\begin{defi}
	Let $R \in \{\OO,k\}$, and let $A$ be any $R$-algebra that is free and finitely generated as an $R$-module
	(i. e. an order if $R=\OO$). Then define the \emph{Picard group} of $A$ as follows:
	\begin{equation}
		\Pic_R(A) := \{ \textrm{ Isomorphism classes of invertible $A^{\op}\otimes_R A$-modules } \}
	\end{equation}
	Note that we will always identify  $A^{\op}\otimes _R B$-modules with the  corresponding $A$-$B$-bimodules. An $A^{\op}\otimes_R A$-module $X$ 
	is invertible (by definition) if it is projective as a left and as a right $A$-module and there is an $A^{\op}\otimes_R A$-module $Y$ (also projective as a left and as a right $A$-module) such that $X\otimes_A Y \iso Y\otimes_A X \iso {_AA_A}$. Now $\Pic_R(A)$ becomes a group with
	``$-\otimes_A =$'' as its product. 

    Similarly define the \emph{derived Picard group} of $A$ as follows:
    \begin{equation}
	   {\rm TrPic}_R(A) := \{  \textrm{ Isomorphism classes of invertible objects in $\mathcal D^b(A^{\op}\otimes_R A)$ } \}
    \end{equation}
    We say a complex $X$ in $\mathcal{D}^b(A^{\op}\otimes_R A)$ is invertible if there is a $Y$ in 
    $\mathcal{D}^b(A^{\op}\otimes_R A)$ such that $X\otimes^{\mathbb L}_A Y \iso Y \otimes^{\mathbb L}_A X \iso 0\rightarrow {_AA_A} \rightarrow 0$.
    ${\rm TrPic}_R(A)$ is a group with ``$-\otimes^{\mathbb L}_A =$'' as its product.
\end{defi}

\begin{remark}
	Situation as above. There is a group homomorphism
	\begin{equation}
		(\Aut_R(A), \circ) \rightarrow (\Pic_R(A), \otimes_A): \ \alpha \mapsto {_{\id} A _{\alpha}} \iso {_{\alpha^{-1}} A _{\id}}
	\end{equation}
	The kernel of this homomorphism consists of all inner automorphisms of $A$, and we denote its image by $\Out_R(A)$. 
	In case $R=k$ is an algebraically closed field, 
	$\Out_k(A)$ is a linear algebraic group defined over $k$, and we denote its 
	connected component by $\Out^0_k(A)$.
\end{remark}

\begin{remark}
	Situation as above. If $X \in \Pic_R(A)$, then $X$ is projective and finitely generated as a left $A$-module and
	as a right $A$-module. Hence if $P\in\projC_A$, then $P\otimes_A X$ is again in $\projC_A$. If 
	$P$ is indecomposable, then so is $P\otimes_A X$, since $X$ is invertible. This implies that
	there is a group homomorphism from the Picard group of $A$ into the symmetric group ${\rm Sym}(\mathcal P)$ on $\mathcal P$
	\begin{equation}
		\Pic_R(A) \longrightarrow {\rm Sym}(\mathcal{P}): \ X \mapsto \left[P \mapsto P \otimes_A X\right]
	\end{equation}
	where $\mathcal{P}$ is the set of all isomorphism classes of finitely generated projective indecomposable 
	$A$-modules. Define $\Pic^s_R(A)$ to be the kernel of this group homomorphism, and $\Out^s_R(A)$ to be the intersection of $\Pic^s_R(A)$ with $\Out_R(A)$. Define $\Aut_R^s(A)$ to be the preimage of 
	$\Out_R^s(A)$ under the canonical epimorphism $\Aut_R(A) \twoheadrightarrow \Out_R(A)$.
\end{remark}

\begin{remark}
    Situation as above.
    Then we get a series of embeddings
    \begin{equation}
	   \Out_R(A) \hookrightarrow \Pic_R(A) \hookrightarrow {\rm TrPic}_R(A)
    \end{equation}
\end{remark}

\begin{notation}
	If $A$ is a ring and $T\in \mathcal K^b(\projC_A)$ is a tilting complex with endomorphism ring $B$, then we denote 
	by
	\begin{equation}
		\mathcal G_T: \ \mathcal D^b(A) \stackrel{\sim}{\longrightarrow} D^b(B)
	\end{equation}
	an equivalence which agrees on objects with taking $T$-resolutions. For proper definitions and 
	proof of the existence of such an equivalence we refer the reader to \cite{RickardOne}. For our purposes, the most important property of
	$\mathcal G_T$ is that $\mathcal G_T(T)$ is isomorphic to the stalk complex $0\rightarrow B \rightarrow 0$.
\end{notation}

\begin{remark}
	Let $\Lambda$ be an $\OO$-order.
	The functor $k\otimes_\OO -: \modC_{\Lambda} \longrightarrow \modC_{k\otimes\Lambda}$ has a (unique) left-derived functor
	$k\otimes^{\mathbb L}_\OO -: \mathcal D^-(\Lambda) \longrightarrow \mathcal D^-(k\otimes\Lambda)$, which restricts to a functor from $\mathcal K^b(\projC_\Lambda)$ to $\mathcal K^b(\projC_{k\otimes \Lambda})$. For a complex  $C\in\mathcal K^b(\projC_\Lambda)$,
	 $k\otimes^{\mathbb L}_\OO C$ is obtained by simply 
	 applying
	 $k\otimes_\OO -$ to this complex viewed as a sequence of modules. Hence, for objects $C\in\mathcal K^b(\projC_\Lambda)$, there
	 is no harm in writing $k\otimes_\OO C$ instead of $k\otimes_\OO^{\mathbb L} C$.
\end{remark}

\begin{remark}\label{remark_tilting_stalk2stalk}
	Let $A$ and $B$ be $R$-algebras and let $\mathcal{F}: \ \mathcal{D}^b(A)\longrightarrow \mathcal D^b(B)$ be an
	equivalence that sends the stalk complex $0\rightarrow A \rightarrow 0$ to $0\rightarrow B \rightarrow 0$. Then there is an $\alpha: \ A \stackrel{\sim}{\longrightarrow} B$ such that
	$\mathcal{F}(X) \iso X\otimes^{\mathbb L}_A {_\alpha B _{\id}}$ for all objects $X \in \mathcal D^b(A)$.
	This follows from \cite[Proposition 7.1]{RickardOne}.
\end{remark}

\begin{lemma}\label{lemma_twosided_tilt_ex}
	Let $A$ be a finite-dimensional $k$-algebra and $T \in \mathcal{K}^b(A)$ a tilting complex
	with endomorphism ring $B$. Then there
	exists a two-sided tilting complex $_B X_A \in \mathcal D^b(B^{\op}\otimes_k A)$ with restriction to
	$\mathcal{D}^b(A)$ isomorphic to $T$. 
\begin{proof}
	By \cite{RickardOne}, there exists a functor $\mathcal{F}: \ \mathcal D^b(B) \longrightarrow \mathcal D^b(A)$ sending $0\rightarrow B \rightarrow 0$ to $T$. By \cite[Corollary 3.5]{RickardDerEqDerFun}
	this equivalence is afforded by ${\rm RHom}_B(Y,-)$ for some $Y\in \mathcal D^b({A^{\op}}\otimes B)$.
	This $Y$ has an inverse $X\in \mathcal{D}^b({B^{\op}}\otimes A)$ such that ${\rm RHom}_B(Y,-)\iso -\otimes_B^{\mathbb{L}}X$ (see \cite[Definition 4.2]{RickardDerEqDerFun} and the remarks following it). Since $B \otimes^{\mathbb L}_B X \iso \mathcal{F}(B)  \iso T$, $X$ has the desired properties.
\end{proof}
\end{lemma}

\begin{lemma}\label{lemma_inverse_two_term}
	Let $A$ be a finite-dimensional symmetric  $k$-algebra and let
	\begin{equation}
		T = 0 \rightarrow P_1 \rightarrow P_0 \rightarrow 0
	\end{equation}
	be a two-term tilting complex. Then $\mathcal{G}_{T}(0 \rightarrow A \rightarrow 0)$  is again
	a two-term tilting complex.
\begin{proof}
	Set $B := \End_{\mathcal D^b(A)}(T)$. Let $X \in \mathcal D^b(B^{\op}\otimes A)$ be a two-sided tilting complex with restriction to $\mathcal D^b(A)$ isomorphic to $T$. Let $Y\in\mathcal D^b(A^{\op}\otimes B)$ be the inverse of $X$.
	By \cite[Lemma 9.2.6]{DerEq} we may assume that $X$ is a bounded complex of $A$-$B$-bimodules that 
	become projective upon restriction to $A$ and restriction to $B$. We may then furthermore assume that 
	$Y = \Hom_k(X, k)$ (see \cite[Corollary 9.2.5]{DerEq}; Note that
	both Lemma 9.2.6 and Corollary 9.2.5 in \cite{DerEq} use that $A$ is symmetric). Hence $Y$
	has non-vanishing homology in precisely two adjacent degrees, since the same can be said about $X$ and $\Hom_k(-,k)$ is exact on vector spaces. Now $-\otimes^{\mathbb{L}}_A Y$
	sends $T$ to $0\rightarrow B\rightarrow 0$, which implies that for some automorphism $\gamma: B \rightarrow B$ the functor $-\otimes^{\mathbb{L}}_A Y \otimes ^{\mathbb L}_B {_{\id}B_{\gamma}}$ agrees with $\mathcal G_{T}(-)$ on objects. Hence the image of $0 \rightarrow A \rightarrow 0$
	under $\mathcal G_{T}(-)$ is equal to the restriction $ Y \otimes ^{\mathbb L}_B {_{\id}B_{\gamma}}$ to $\mathcal D^b(B)$. Therefore it is a bounded complex of projective $B$-modules that has non-zero homology (at most) in two (adjacent) degrees. Since $A$ and $B$ are symmetric (so in particular self-injective), any injection of a projective module and any epimorphism onto a projective module splits. Hence $ Y \otimes ^{\mathbb L}_B {_{\id}B_{\gamma}}$ is isomorphic in
	$\mathcal{K}^b(\projC_B)$ to a two-term complex.
\end{proof}
\end{lemma}

For the rest of the section, let $k$ be algebraically closed.

\begin{thm}[Rouquier, Huisgen-Zimmermann, Saor\'{\i}n]\label{thm_out0_inv}
	Let $A, B$ be finite-dimensional $k$-algebras and $X$ a bounded complex of $A$-$B$-bimodules inducing
	an equivalence between $\mathcal D^b(A)$ and $\mathcal D^b(B)$ (i. e., a two-sided tilting complex). 
	Then there exists a (unique) isomorphism of algebraic groups
	\begin{equation}
		\sigma:\ \Out_k^0(A) \stackrel{\sim}{\longrightarrow} \Out_k^0(B) 
	\end{equation}
	such that
	\begin{equation}
		_{\id} A _{\alpha} \otimes^{\mathbb L}_A X \iso X \otimes^{\mathbb L}_B {_{\id} B _{\sigma(\alpha)}}
	\end{equation}
	for all $\alpha \in \Out_k^0(A)$.
\begin{proof}
	The Theorem was stated in this form in \cite[Theorem 3.4]{RouquierICM2006}. A proof can be found in
	\cite{ZimmermannGeometryChainComplexes} or in \cite{RouquierAutomorphFrench}.
\end{proof}
\end{thm}

\begin{thm}[Jensen, Su, Zimmermann]\label{thm_jensen_su_z}
	Let $A$ be a finite-dimensional $k$-algebra. Then up to isomorphism in $\mathcal{K}^b(\projC_A)$ there
	exists at most one two-term (partial) tilting complex
	\begin{equation}
		0 \rightarrow P_1 \rightarrow P_0 \rightarrow 0
	\end{equation}
	with fixed homogeneous components $P_0$ and $P_1$.
\begin{proof}
	See \cite[Corollary 8]{JensenXuZDegenerations}.
\end{proof}
\end{thm}

\begin{corollary}
	Let $A$ be a finite-dimensional $k$-algebra and $T$ a tilting complex over $A$. Then
	\begin{enumerate}
	 \item $T\otimes_A {_\id A_{\gamma}}\iso T$ for all $\gamma\in\Out_k^0(A)$.
	 \item If $T$ is a two-term complex, then $T\otimes_A {_\id A_{\gamma}}\iso T$
		for all $\gamma\in\Out_k^s(A)$.
	\end{enumerate}
\begin{proof}
	The first point follows from Theorem \ref{thm_out0_inv} and Lemma \ref{lemma_twosided_tilt_ex}.
	The second point follows from Theorem \ref{thm_jensen_su_z} and the definition of 
	$\Out_k^s(A)$.
\end{proof}
\end{corollary}

\section{A Correspondence of Lifts}\label{section_corr_lifts}

In this section we introduce a bijection between ``lifts'' of derived equivalent finite-dimensional $k$-algebras.

\begin{defi}\label{defi_lifts}
	For a finite-dimensional $k$-algebra $\overline\Lambda$ define its \emph{set of lifts} as follows:
	\begin{equation}
		\Lifts(\overline{\Lambda}) := \left\{ (\Lambda, \varphi) \ | \  \textrm{$\Lambda$ is an $\OO$-order and $\varphi: k\otimes\Lambda \stackrel{\sim}{\rightarrow} \overline \Lambda$ is an isomorphism}\right\}
		\bigg / \sim
	\end{equation}
	where we say $(\Lambda, \varphi) \sim (\Lambda',\varphi')$ if and only if 
	\begin{enumerate}
	 \item There is an isomorphism $\alpha: \Lambda \stackrel{\sim}{\rightarrow} \Lambda'$
	 \item There is a $\beta \in \Aut_k(\overline\Lambda)$ such that the functor  $-\otimes^{\mathbb L}_{\overline\Lambda} {_\beta \overline\Lambda_\id}$ fixes all isomorphism classes of tilting complexes in $\mathcal K^b(\projC_{\overline\Lambda})$
	\end{enumerate}
	such that $\varphi=\beta\circ\varphi'\circ(\id_k\otimes \alpha)$.
\end{defi}
Our bijection will be based on the following theorem of Rickard:
\begin{thm}[{\cite[Theorem 3.3.]{RickardLiftTilting}}]\label{thm_unique_lift_tilt}
	Let $\Lambda$ be an $\OO$-order and let  $\overline{T} \in \mathcal{K}^b(\projC_{k\otimes \Lambda})$ be a tilting complex for
	$k\otimes \Lambda$. Then there exists a unique (up to isomorphism in $\mathcal{D}^b(\Lambda)$) tilting
	complex $T \in \mathcal{K}^b(\projC_\Lambda)$ with $k\otimes T \iso \overline{T}$.
	$\End_{\mathcal{D}^b(\Lambda)}(T)$ is torsion-free and
	\begin{equation}
		k\otimes \End_{\mathcal{D}^b(\Lambda)}(T) \iso  \End_{\mathcal{D}^b(k\otimes \Lambda)}(\overline{T}) 
	\end{equation}
\end{thm}

\begin{remark}\label{remark_partial_tilting}
	By \cite[Proposition 3.1.]{RickardLiftTilting} it is immediately clear that we can replace the word ``tilting complex'' by ``partial tilting complex'' in the above theorem (where we understand ``partial tilting complex'' as defined in \cite[Definition 3.2.1.]{DerEq}). 
\end{remark}

\begin{lemma}\label{lemma_commuting_reduction_modp}
    Let $\Lambda$ be an $\OO$-order and $T\in\mathcal K^b(\projC_\Lambda)$ a tilting complex. Define $\Gamma := \End_{\mathcal D^b(\Lambda)}(T)$, and assume that $\Gamma$ is also an $\OO$-order. Then $k\otimes T$ is a tilting complex and the $k$-algebras $k\otimes \Gamma$ and $\End_{\mathcal D^b(\Lambda)}(k\otimes T)$
    are (canonically) isomorphic. Moreover, the diagram
    \begin{equation}
	   \xymatrix{
		  \mathcal D^-(\Lambda) \ar[rr]^{\mathcal G_T}  \ar[d]^{k\otimes^{\mathbb L}-} & & \mathcal D^-(\Gamma) \ar[d]^{k\otimes^{\mathbb L}-} \\
		  \mathcal D^-(k\otimes\Lambda) \ar[rr]^{\mathcal G_{k\otimes T}} & & \mathcal D^-(k\otimes\Gamma)
	   }
    \end{equation}
    commutes on objects.
\begin{proof}
  This follows from {\cite[Proposition 2.4.]{RickardDerEqDerFun}}.
\end{proof}
\end{lemma}

For the rest of the section let $\overline\Lambda$ and $\overline\Gamma$ be two derived equivalent finite-dimensional $k$-algebras. Furthermore let $X\in \mathcal{D}^b(\overline\Lambda^{\op} \otimes_k \overline\Gamma)$ be a two-sided tilting complex, and let $X^{-1}$ be its inverse. Let $\overline T$ be the restriction of $X^{-1}$ to $\mathcal D^b(\projC_{\overline\Lambda})$ and likewise let $\overline S$ be the restriction of $X$ to $\mathcal D^b(\projC_{\overline\Gamma})$.

\begin{defi}\label{defi_phi}
 	Define a map
	\begin{equation}
		\Phi_X: \ \Lifts(\overline \Lambda) \longrightarrow \Lifts(\overline \Gamma) 
	\end{equation}
	as follows: Assume $(\Lambda,\varphi)\in \Lifts(\overline\Lambda)$.
	Let $T$ be the lift of $\overline T \otimes_{\overline\Lambda} { _\id \overline \Lambda _\varphi}$ (which exists and is unique 
	 by Theorem \ref{thm_unique_lift_tilt}). 
	We put $\Phi_X(\Lambda,\varphi) = (\Gamma, \psi)$, where
	$\Gamma = \End_{\mathcal D^b(\Lambda)}(T)$ and $\psi: k\otimes \Gamma \stackrel{\sim}{\rightarrow} \overline\Gamma$ is an isomorphism
	 such that the following diagram commutes on objects:
	 \begin{equation}\label{eqn_def_PHI}
	 \xymatrix{
	  \mathcal{D}^-(\Lambda) \ar[rr]^{\mathcal{G}_{T}(-)}\ar[d]_{k\otimes^{\mathbb{L}}-} &  & \mathcal{D}^-(\End_{\mathcal D^b(\Lambda)}(T)) 
	 \ar[d]^{k\otimes^{\mathbb{L}}-} \ar@{=}[rr]& & \mathcal{D}^-(\Gamma)\ar[d]^{k\otimes^{\mathbb{L}}-} \\
	 \mathcal{D}^-({k\otimes\Lambda}) \ar[d]_{-\otimes^{\mathbb{L}}_{k\otimes\Lambda} {_\varphi \overline{\Lambda}_\id}}
		\ar[rr]^(.4){\mathcal{G}_{k\otimes T}(-)} & & \mathcal{D}^-(\End_{\mathcal D^b(k\otimes\Lambda)}(k\otimes T))\ar[d]^{\mathcal E} 
		\ar@{=}[rr]
		& & \mathcal{D}^-(k\otimes\Gamma)\ar[d]^{-\otimes^{\mathbb L}_{k \otimes \Gamma} {_{\psi} \overline\Gamma_{\id}}} \\
	 \mathcal{D}^-({\overline\Lambda})\ar[rr]_{-\otimes^{\mathbb L}_{\overline\Lambda} X} & & \mathcal{D}^-({\overline\Gamma}) \ar@{=}[rr] & & \mathcal{D}^-(\overline{\Gamma}) \\
    	}		
	\end{equation}
	 Here, $\mathcal E$ is defined so that the bottom left square commutes.
\begin{proof}[Proof of well-definedness]
	 First note that the top left square commutes on objects by Lemma \ref{lemma_commuting_reduction_modp}. Thus the left half
	 of the diagram will commute on objects. Note furthermore that $\mathcal E$ sends $0\rightarrow k\otimes\Gamma \rightarrow 0$
	 to $0\rightarrow \overline\Gamma \rightarrow 0$, and hence a $\psi$ making the diagram commutative on objects
	 can be chosen due to Remark \ref{remark_tilting_stalk2stalk}. This $\psi$ is unique up to an automorphism $\beta$
	 of $\overline\Gamma$ such that $-\otimes^{\mathbb L}_{\overline\Gamma}{_\id}\overline\Gamma_\beta$ fixes all objects
	 in $\mathcal D^-(\overline\Gamma)$, and hence in particular fixes all tilting complexes. Therefore
	 the equivalence class of $(\Gamma,\psi)$ is certainly independent of the particular choice of $\psi$.

	 Now assume $(\Lambda,\varphi)\sim (\Lambda', \varphi')$, that is, there are $\alpha$ and $\beta$ as in Definition \ref{defi_lifts}
	 such that $\varphi = \beta \circ \varphi' \circ (\id_k\otimes\alpha)$.  We need to show that $(\Gamma,\psi) :=\Phi_X(\Lambda, \varphi)\sim\Phi_X(\Lambda',\varphi')=: (\Gamma',\psi')$, $\Phi_X$ being given by the construction above. We get the following diagram (where we define $T'$ analogous to $T$):
	 \begin{equation}
	   \xymatrix{
	   \mathcal{D}^-(\Gamma') \ar[d]^{k\otimes^{\mathbb L}-} & & \ar[d]^{k\otimes^{\mathbb L}-} \ar[ll]_{\mathcal G_{T'}(-)}\mathcal{D}^-(\Lambda')\ar[rr]^{-\otimes^{\mathbb L}_{\Lambda'} {_\id \Lambda'_\alpha} } & & \ar[d]^{k\otimes^{\mathbb L}-} \mathcal{D}^-(\Lambda) \ar[rr]^{\mathcal G_T(-)} & & \ar[d]^{k\otimes^{\mathbb L}-}\mathcal{D}^-(\Gamma) \\
	   \mathcal{D}^-(k\otimes \Gamma') \ar[d]^{-\otimes^{\mathbb L}_{k\otimes\Gamma'} {_{\psi'}\overline\Gamma_\id}}& & \mathcal{D}^-(k\otimes\Lambda') \ar[ll]_{\mathcal G_{k\otimes T'}(-)}  \ar[rr]_{-\otimes^{\mathbb L}_{k\otimes \Lambda'} {_\id k\otimes\Lambda'_{\id_k\otimes \alpha}}}\ar[d]^{-\otimes^{\mathbb L}_{k\otimes\Lambda'} {_{\varphi'}\overline\Lambda_\id}} & & \ar[d]^{-\otimes^{\mathbb L}_{k\otimes\Lambda} {_{\varphi}\overline\Lambda_\id}} \mathcal{D}^-(k\otimes\Lambda) \ar[rr]^{\mathcal G_{k\otimes T}(-)} & & \mathcal{D}^-(k\otimes\Gamma)\ar[d]^{-\otimes^{\mathbb L}_{k\otimes\Gamma} {_{\psi}\overline\Gamma_\id}} \\
	   \mathcal{D}^-(\overline\Gamma) & & \mathcal{D}^-(\overline\Lambda) \ar@{=}[rr] \ar[ll]^{-\otimes^{\mathbb L}_{\overline\Lambda} X} & & \mathcal{D}^-(\overline\Lambda) \ar[rr]_{-\otimes^{\mathbb L}_{\overline\Lambda} X} & & \mathcal{D}^-(\overline\Gamma) \\
	   }
	 \end{equation}
	 This diagram will commute at the very least on tilting complexes (that is, if we take a tilting complex in any of those
	 categories and take its image under a series of arrows in the above diagram, the isomorphism class of the outcome will not depend on the path we have chosen). 
	 Note that all horizontal arrows are equivalences, and so we get a diagram (again commutative on tilting complexes)
	 \begin{equation}
		\xymatrix{
		  \mathcal{D}^-(\Gamma') \ar[d]^{k\otimes^{\mathbb L}-} \ar[rr]^{\mathcal F_1(-)}& & \mathcal{D}^-(\Gamma) \ar[d]^{k\otimes^{\mathbb L}-} \\ 
		  \mathcal{D}^-(k\otimes\Gamma')\ar[d]^{-\otimes^{\mathbb L}_{k\otimes\Gamma'} {_{\psi'}\overline\Gamma_\id}} \ar[rr]^{\mathcal F_2(-)} & & \mathcal{D}^-(k\otimes\Gamma) \ar[d]^{-\otimes^{\mathbb L}_{k\otimes\Gamma} {_{\psi}\overline\Gamma_\id}} \\
		  \mathcal{D}^-(\overline\Gamma) \ar@{=}[rr] & & \mathcal{D}^-(\overline\Gamma)  \\
		}
	 \end{equation}
	 where $\mathcal F_1$ and $\mathcal F_2$ are two equivalences.
	 Due to commutativity on tilting complexes, $\mathcal F_2$ needs to send $0\rightarrow k\otimes \Gamma' \rightarrow 0$ to
	 $0\rightarrow k\otimes \Gamma \rightarrow 0$. Due to unique lifting (and again commutativity), $\mathcal F_1$ needs to 
	 send $0\rightarrow \Gamma' \rightarrow 0$ to $0\rightarrow \Gamma \rightarrow 0$. Hence there is an isomorphism
	 $\alpha: \Gamma' \rightarrow \Gamma$ such that $\mathcal F_1(-)$ agrees on objects with 
	 $-\otimes^{\mathbb L}_{\Gamma'} {_\alpha \Gamma _\id}$. Due to commutativity, $\mathcal F_2(-)$ then needs to agree on tilting complexes with $-\otimes^{\mathbb L}_{k\otimes \Gamma'} {_{\id_k\otimes \alpha} k\otimes \Gamma_\id}$ (this is owed to the fact that every tilting complex lies in the image of $k\otimes^{\mathbb L}-$ due to Theorem \ref{thm_unique_lift_tilt}). Commutativity on 
	 tilting complexes of the lower square then implies that $\psi' = \beta\circ\psi\circ(\id_k\otimes \alpha)$ for some $\beta\in \Aut_k(\overline \Gamma)$ so that
	 $-\otimes^{\mathbb L}_{\overline\Gamma}{_\beta \overline\Gamma_\id}$ fixes all tilting complexes in $\mathcal K^b(\projC_{\overline \Gamma})$.
	 By definition this means $(\Gamma,\psi)\sim(\Gamma',\psi')$.
\end{proof}
\end{defi}

\begin{prop}
	The maps $\Phi_X$ and $\Phi_{X^{-1}}$ are mutually inverse. In particular, they induce a bijection
	\begin{equation}
		\Lifts(\overline\Lambda) \longleftrightarrow \Lifts(\overline\Gamma)
	\end{equation}
\begin{proof}
	 We keep the notation of Definition \ref{defi_phi}. Set $(\Gamma,\psi) := \Phi_X(\Lambda,\varphi)$ and
	 $(\Lambda',\tilde\varphi) := \Phi_{X^{-1}}(\Gamma,\psi)$.
	 Furthermore, let $S$ be the lift of $\overline S\otimes_{\overline\Gamma} {_\id \overline\Gamma_\psi}$.
	Consider the following diagram (obtained by composing (\ref{eqn_def_PHI}) with itself)
	\begin{equation}
	 \xymatrix{
	  \mathcal{D}^-(\Lambda) \ar[rr]^{\mathcal{G}_{T}(-)}\ar[d]_{k\otimes^{\mathbb{L}}-} &  & \mathcal{D}^-(\Gamma) 
	 \ar[d]^{k\otimes^{\mathbb{L}}-} \ar[rr]^{\mathcal{G}_S(-)}& & \mathcal{D}^-(\Lambda')\ar[d]^{k\otimes^{\mathbb{L}}-} \\
	 \mathcal{D}^-({k\otimes\Lambda}) \ar[d]_{-\otimes^{\mathbb{L}}_{k\otimes\Lambda} {_\varphi \overline{\Lambda}_\id}}
		\ar[rr]^{\mathcal{G}_{k\otimes T}(-)} & & \mathcal{D}^-(k\otimes\Gamma)\ar[d]^{\otimes^{\mathbb L}_{k\otimes\Gamma}{_\psi\overline\Gamma_\id}} 	 
		\ar[rr]^{\mathcal G_{k\otimes S}(-)}
		& & \mathcal{D}^-(k\otimes\Lambda')\ar[d]^{-\otimes^{\mathbb L}_{k \otimes \Lambda'} {_{\tilde\varphi} \overline\Lambda_{\id}}} \\
	 \mathcal{D}^-({\overline\Lambda})\ar[rr]_{-\otimes^{\mathbb L}_{\overline\Lambda} X} & & \mathcal{D}^-({\overline\Gamma}) \ar[rr]_{-\otimes^{\mathbb L}_{\overline\Gamma} X^{-1}} & & \mathcal{D}^-(\overline{\Gamma}) \\
    	}		
	\end{equation}
	Commutativity on objects implies that $\mathcal G_{k\otimes S} \circ \mathcal G_{k\otimes T}$
	sends the stalk complex $0\rightarrow k \otimes \Lambda\rightarrow 0$ to the stalk complex $0\rightarrow k\otimes \Lambda'\rightarrow 0$. 
	 Hence, due to unique lifting, $\mathcal{G}_S \circ \mathcal G_T$ sends $0\rightarrow \Lambda \rightarrow 0$ to $0\rightarrow \Lambda' \rightarrow 0$. $\mathcal{G}_S \circ \mathcal G_T$ hence agrees
	on objects with $-\otimes^{\mathbb L}_{\Lambda} {_\id \Lambda _\alpha}$ for some isomorphism $\alpha: \Lambda'\stackrel{\sim}{\longrightarrow} \Lambda$.
	Thus, $\mathcal G_{k\otimes S} \circ \mathcal G_{k\otimes T}$ will agree on tilting complexes with
	$-\otimes^{\mathbb L}_{k\otimes\Lambda} {_\id k\otimes\Lambda_{\id_k \otimes \alpha}}$. Hence, due to commutativity on objects, we must have that $-\otimes^{\mathbb L}_{k\otimes\Lambda} {_{\tilde\varphi} \overline\Lambda _\id}$
	and $-\otimes^{\mathbb L}_{k\otimes\Lambda} {_{\id_k\otimes\alpha} k\otimes\Lambda _\id} \otimes^{\mathbb L}_{k\otimes\Lambda} {_\varphi \overline\Lambda_\id} = - \otimes^{\mathbb L}_{k\otimes\Lambda} {_{\varphi\circ (\id_k\otimes \alpha)} \overline\Lambda_\id}$ agree on tilting complexes. This however is the same as saying that 
	$\tilde \varphi = \beta \circ \varphi \circ (\id_k\otimes \alpha)$, where $\beta\in\Aut_k(\overline\Lambda)$ is an automorphism such that $-\otimes^{\mathbb L}_{\overline \Lambda}{_{\beta}\overline\Lambda_\id}$ fixes all tilting complexes. 
     This means, by definition, that $(\Lambda, \varphi)\sim (\Lambda', \tilde\varphi)$. So we proved that $\Phi_{X^{-1}}\circ\Phi_X=\id$,
	and $\Phi_{X}\circ\Phi_{X^{-1}}=\id$ follows by swapping the roles of $X$ and $X^{-1}$.
\end{proof}
\end{prop}

\begin{prop}
	$\Out_k(\overline\Lambda)$ acts on $\Lifts(\overline\Lambda)$ from the left via
	\begin{equation}
		\alpha\cdot (\Lambda, \varphi) := (\Lambda, \alpha \circ \varphi)
	\end{equation}
\begin{proof}
 	The above formula clearly defines an action of $\Aut_k(\overline\Lambda)$. In order to verify that it defines an action of $\Out_k(\overline{\Lambda})$, we just need to check that
	for any inner automorphism $\alpha$ of $\overline{\Lambda}$ we have
	$(\Lambda,\varphi)\sim (\Lambda,\alpha\circ\varphi)$. But an inner automorphism  $\alpha$ of $\overline\Lambda$
	gives us an inner automorphism $\varphi^{-1}\circ\alpha\circ\varphi$ of $k\otimes\Lambda$, which lifts to an inner automorphism $\hat\alpha$ of $\Lambda$ 
	(since the natural map of unit groups $\Lambda^\times \rightarrow (k\otimes \Lambda)^\times$ is surjective).
	$(\Lambda,\varphi)$ and $ (\Lambda,\alpha\circ\varphi)=
    (\Lambda,\varphi\circ(\id_k\otimes\hat \alpha))$ are then clearly equivalent in the sense of Definition \ref{defi_lifts}.
\end{proof}
\end{prop}

\begin{prop}\label{prop_out0_trivial_on_lifts}
    If $k$ is algebraically closed, then $\Out_k^0(\overline\Lambda)$ lies in the kernel of the action of 
    $\Out_k(\overline\Lambda)$ on $\Lifts(\overline\Lambda)$.
\begin{proof}
    This follows directly from Theorem \ref{thm_out0_inv}.
\end{proof}
\end{prop}

\begin{prop}
	 Let $\Out_k(\overline{\Lambda})_{\overline T}$ respectively $\Out_k(\overline{\Gamma})_{\overline S}$
	 denote the stabilizers of the isomorphism classes of $\overline T$ respectively $\overline S$.
	There is an isomorphism 
	\begin{equation}
		-^X: \Out_k(\overline{\Lambda})_{\overline T} \stackrel{\sim}{\longrightarrow} \Out_k(\overline{\Gamma})_{\overline S}
	\end{equation}
	such that for all $\alpha\in \Out_k(\overline{\Lambda})_{\overline T}$ we have
	\begin{equation}\label{eqn_hdjhduio}
		\Phi_X(\alpha\cdot(\Lambda,\varphi))=\alpha^X\cdot\Phi_X(\Lambda,\varphi)
	\end{equation}
	In particular, $\Phi_X$ induces a bijection
	\begin{equation}\label{eqn_jkj373}
		\Out_k(\overline{\Lambda})_{\overline T} \setminus \Lifts(\overline\Lambda) \longleftrightarrow \Out_k(\overline{\Gamma})_{\overline S} \setminus \Lifts(\overline\Gamma)
	\end{equation}
\begin{proof}
	 Set
  \begin{equation}
	 -^X: \Out_k(\overline\Lambda)_{\overline T} \longrightarrow {\rm TrPic}(\overline \Gamma): \alpha \mapsto X^{-1} \otimes^{\mathbb L}_{\overline\Lambda} {_\id\overline\Lambda_\alpha} \otimes^{\mathbb L}_{\overline\Lambda} X
  \end{equation}
  First note that the restriction of $X^{-1}$ to $\mathcal D^b(\Lambda)$ is isomorphic to $\overline T$ by definition
  of $\overline T$. Since $\alpha$ stabilizes the isomorphism class of $\overline T$, the restriction of $X^{-1} \otimes^{\mathbb L}_{\overline\Lambda} {_\id\overline\Lambda_\alpha} \otimes^{\mathbb L}_{\overline\Lambda} X$ to $\mathcal D^b(\overline\Gamma)$ is isomorphic to $0 \rightarrow \overline \Gamma \rightarrow 0$. Thus $X^{-1} \otimes^{\mathbb L}_{\overline\Lambda} {_\id\overline\Lambda_\alpha} \otimes^{\mathbb L}_{\overline\Lambda} X$ is isomorphic to $0\rightarrow{_\id \overline\Gamma_\beta}\rightarrow 0$ for some $\beta \in\Aut_k(\overline\Gamma)$. That is, the image of $-^X$ as defined above is indeed contained in $\Out_k(\overline\Gamma) \leq {\rm TrPic}(\overline \Gamma)$. Now $\overline S$ is by definition just the restriction of $X$ to $\mathcal D^b(\overline\Gamma)$, and hence
  $\overline S\otimes^{\mathbb L}_{\overline\Gamma}X^{-1} \otimes^{\mathbb L}_{\overline\Lambda} {_\id\overline\Lambda_\alpha} \otimes^{\mathbb L}_{\overline\Lambda} X$ is isomorphic to the restriction of ${_\id\overline\Lambda_\alpha} \otimes^{\mathbb L}_{\overline\Lambda} X$ to $\mathcal D^b(\overline\Gamma)$
  which is again isomorphic to $\overline S$. So $-^X$ does indeed define a map with image contained in $\Out_k(\overline\Gamma)_{\overline S}$.
  It is also easy to see that $-^X$ is a group homomorphism, and that $-^{X^{-1}}$ is a two-sided inverse for $-^X$.

  Now the claim of (\ref{eqn_hdjhduio}) follows from the commutativity of the following diagram
  \begin{equation}
	 \xymatrix{
		\mathcal D^-(\overline \Lambda) \ar[d]_{-\otimes^{\mathbb L}_{\overline \Lambda} {_\id\overline\Lambda_{\alpha}}} \ar[rr]^{-\otimes^{\mathbb L}_{\overline \Lambda} X} && \mathcal D^-(\overline \Gamma) \ar[d]^{-\otimes^{\mathbb L}_{\overline \Gamma} {_\id\overline\Gamma_{\alpha^X}}} \\	 
		\mathcal D^-(\overline \Lambda) \ar[rr]_{-\otimes^{\mathbb L}_{\overline \Lambda} X} && \mathcal D^-(\overline \Gamma) 
	 }
  \end{equation}
  by gluing it below diagram (\ref{eqn_def_PHI}).
\end{proof}
\end{prop}


\begin{defi}
    Define the set $\mathfrak L(\overline \Lambda)$ to be the set of all isomorphism classes of $\OO$-orders $\Lambda$ such that
    $k\otimes \Lambda \iso \overline\Lambda$.  Clearly, $\mathfrak L(\overline \Lambda)$ is in bijection with $\Out_k(\overline\Lambda) \setminus \Lifts(\overline\Lambda)$. Furthermore, we define the projection map
    \begin{equation}
	   \Pi: \Lifts(\overline\Lambda) \longrightarrow \mathfrak L(\overline \Lambda): (\Lambda,\varphi) \mapsto \Lambda
    \end{equation}
\end{defi}

\begin{corollary}\label{corr_lifts_two_term}
    Assume $k$ is algebraically closed, $\overline\Lambda$ is symmetric, and $\overline T$ is a two-term complex. 
    Assume furthermore that $\Out_k^s(\overline\Lambda)=\Out_k(\overline\Lambda)$ and $\Out_k^s(\overline\Gamma)=\Out_k(\overline\Gamma)$
    (a sufficient criterion for this is for instance that the Cartan matrices of $\overline\Lambda$ and $\overline\Gamma$ have no
    non-trivial permutation symmetries).
  Then there is a bijection 
    \begin{equation}
	\mathfrak L (\overline\Lambda) \longleftrightarrow \mathfrak L (\overline\Gamma)
    \end{equation}
\begin{proof}
    Lemma \ref{lemma_inverse_two_term} implies that $\overline S$ may be assumed to be a two-term complex as well.
    The assertion now follows from (\ref{eqn_jkj373}) together with Theorem \ref{thm_jensen_su_z}, since the latter implies that
    $\Out_k(\overline\Lambda)_{\overline T}=\Out_k(\overline\Lambda)$ and
    $\Out_k(\overline\Gamma)_{\overline S}=\Out_k(\overline\Gamma)$.
\end{proof}
\end{corollary}

The following proposition is useful to prove a ``unique lifting property'' for the group ring
of $\SL_2(p^f)$ in defining characteristic, which we will do in a later paper.
\begin{prop}
    Assume $k$ is algebraically closed.
    Let $\Lambda \in \mathfrak L(\overline\Lambda)$, and let $\gamma: k\otimes\Lambda\stackrel{\sim}{\rightarrow}\overline\Lambda$.
    be an isomorphism. Now assume
   \begin{equation}\label{eqn_JHKHIOu}
   	\overline{\Aut_\OO(\Lambda)} \cdot \Out_k^0(\overline\Lambda) = \Out_k(\overline\Lambda)
   \end{equation}
    where $\overline{\Aut_\OO(\Lambda)}$ is the image of $\Aut_\OO(\Lambda)$ in $\Out_k(\overline\Lambda)$
    (here we identify $k\otimes\Lambda$ with $\overline\Lambda$ via $\gamma$).
    Then the fiber $\Pi^{-1}(\{\Lambda\})$ has cardinality one.
 \begin{proof}
    Let  $(\Lambda, \varphi)\in\Lifts(\overline\Lambda)$ for some $\varphi: k\otimes\Lambda\stackrel{\sim}{\longrightarrow}\overline\Lambda$.
    Now if (\ref{eqn_JHKHIOu}) holds, we can write $\gamma\circ\varphi^{-1}=\gamma\circ (\id_k\otimes \hat\alpha) \circ \gamma^{-1}\circ \beta$ for some $\hat\alpha\in \Aut_\OO(\Lambda)$ and $\beta \in \Aut_k(\overline\Lambda)$ such that the image of
    $\beta$ in $\Out_k(\overline\Lambda)$ lies in $\Out_k^0(\overline\Lambda)$. Hence $\gamma \circ (\id_k\otimes\hat\alpha^{-1})=  \beta \circ \varphi$.
    Proposition \ref{prop_out0_trivial_on_lifts} (together with the definition of ``$\sim$'') implies $(\Lambda,\gamma)\sim (\Lambda,\beta^{-1}\circ \gamma\circ(\id_k\otimes\hat\alpha^{-1}))=(\Lambda, \varphi)$.
\end{proof}
\end{prop}

\section{Tilting Orders in Semisimple Algebras}

What we would want to do now is to partition the set $\Lifts(\overline \Lambda)$ into manageable pieces, so that the map
$\Phi_X$ defined in the previous section restricts to a bijection between corresponding pieces of $\Lifts(\overline \Lambda)$
and $\Lifts(\overline \Gamma)$. In order to do that in the next section, we first need to study how those properties of orders that
we are interested in behave under derived equivalences.
The results in this section are elementary and therefore certainly known, but we were unable to find explicit references for most of them, which is why we include proofs.
We assume throughout this section that $\Lambda$ is an $\OO$-order in a finite-dimensional semisimple $K$-algebra $A$.

\begin{lemma}\label{lemma_tilting_tensor_K}
	If $T \in \mathcal{K}^b(\projC_\Lambda)$ is a tilting complex for $\Lambda$ and $\End_{\mathcal{D}^b(\Lambda)}(T)$ is
	torsion-free as an $\OO$-module, then $K\otimes T$ is a tilting complex for A. 
	Furthermore, $\End_{\mathcal{D}^b(\Lambda)}(T)$ is a full $\OO$-order in $\End_{\mathcal{D}^b(A)}(K\otimes T)$.
\begin{proof}
	First we show that $\Hom_{\mathcal{D}^b(A)}(K\otimes T, K\otimes T[i]) = 0$ for $i\neq 0$.
	Assume that $\varphi \in \Hom_{\mathcal{D}^b(A)}(K\otimes T, K\otimes T[i])$ for some $i$. Then
	we may view $\varphi$ (or rather a representative of it) as a morphism of graded modules $K\otimes T \rightarrow K\otimes T[i]$
	commuting with the differential. As such we may restrict it to $T$, and for a large enough
	$n\in \N$, we will have $\Im(\pi^n\cdot \varphi) \subseteq T[i]$. Hence $\pi^n \cdot \varphi$
	defines an element in $\Hom_{\mathcal{D}^b(\Lambda)}(T, T[i])$. For $i=0$ this implies
	that $\End_{\mathcal{D}^b(\Lambda)}(T)$ is a full $\OO$-lattice in $\End_{\mathcal{D}^b(A)}(K\otimes T)$. For $i\neq 0$ this implies that $\pi^n\cdot \varphi$ is homotopic to zero, and hence so is 
	$\varphi$ (by dividing the homotopy by $\pi^n$).

	Since $K\otimes^{\mathbb{L}} -: \mathcal{D}^-(\Lambda)\rightarrow \mathcal{D}^-(A)$ is an exact functor between triangulated categories that maps $T$ to $K\otimes T$, it is clear that $\add(K\otimes T)$ contains the image of
	$\add(T)$. But $\add(T)$ is equal to $\mathcal{K}^b(\projC_\Lambda)$ by definition, and so
	in particular contains $0 \rightarrow \Lambda \rightarrow 0$,
	which maps to $0 \rightarrow A \rightarrow 0$, which in turn clearly generates 
	$\mathcal{K}^b(\projC_A)$. Hence $\add(K\otimes T) = \mathcal{K}^b(\projC_A)$.
\end{proof}
\end{lemma}

\begin{lemma}\label{lemma_tilting_wedderburn}
	Situation as above. Let $V_1,\ldots,V_n$ be representatives for the isomorphism classes of simple $A$-modules. Then
	there are sets $\Omega_i$ for $i\in \Z$ with $\biguplus_i \Omega_i = \{1,\ldots, n\}$ and numbers 
	$\delta_j\in \Z_{>0}$ for $j\in\{1,\ldots,n\}$ such that
	\begin{equation}
		K\otimes T \iso_{\mathcal{D}^b(A)} \ldots \stackrel{0}{\rightarrow}
		\underbrace{\bigoplus_{j \in \Omega_i} V_j^{\delta_j}}_{\textrm{degree $i$}} \stackrel{0}{\rightarrow} \bigoplus_{j \in \Omega_{i+1}} V_j^{\delta_j} \stackrel{0}{\rightarrow} \ldots
	\end{equation}
	In particular, each $V_j$ occurs as a direct summand of precisely one of the $H^i(K\otimes T)$.
	Also, it follows that  $H^i(K\otimes T) \iso \bigoplus_{j\in\Omega_i} V_j^{\delta_j}$,
	and the map below is an isomorphism:
	\begin{equation}
		\bigoplus_i H^i:\ \End_{\mathcal{D}^b(A)}(K\otimes T) \stackrel{\sim}{\longrightarrow} \bigoplus_i \End_A(H^i(K\otimes T)) \iso \bigoplus_i \bigoplus_{j\in\Omega_i}\End_A(V_j)^{\delta_j\times \delta_j}
	\end{equation}
\begin{proof}
	$K\otimes T$ is, as a complex over $A$, certainly split, and hence isomorphic in
	the homotopy category to a complex $C$ with differential equal to zero.
	Clearly $H^i(C) = C^i$ and $\End_{\mathcal{D}^b(A)}(C) = \bigoplus_i \End_A(C^i)$.
	So all that remains to show is that any $V_j$ occurs in precisely one $C^i$. But $\Hom_{\mathcal{D}^b(A)}(C,C[l])=0$ for $l\neq 0$ implies that $\Hom_A(C^i,C^{i+l})=0$ for all $l\neq 0$ and hence that $V_j$ occurs in at most one $C^i$. The
	fact that $\add(C)=\mathcal{D}^b(A)$ implies that each $V_j$ has to occur in some $C^i$.
\end{proof}
\end{lemma}

\begin{defi}\label{defi_data_tilting}
	The above lemma contains a definition of sets $\Omega_i$ and numbers $\delta_j$ associated to the tilting complex $T$. We keep this notation. In addition to those, define $\eps: \{1,\ldots,n\}\rightarrow \Z$ to map $j$ to 
	the unique $i$ such that $j \in \Omega_i$. 
   
    Note that in the context of perfect isometries, the numbers $(-1)^{\eps(i)}$ are known as the ``signs'' in the ``bijection with signs''
    induced by a perfect isometry. 
\end{defi}

\begin{thm}[{\cite[Theorem 1]{ZimmTiltedSymmIsSymm}}]\label{thm_tilted_symm_is_symm}
	Assume $\Lambda$ is a symmetric order. Then any $\OO$-algebra $\Gamma$ which is derived equivalent to
	$\Lambda$ is again an $\OO$-order, and symmetric. 
\end{thm}

We close this section by making Theorem \ref{thm_tilted_symm_is_symm} constructive. We wish to give
an explicit symmetrizing form (as defined below) for $\Gamma$, provided we know one for $\Lambda$ 
(which we usually do, for instance in the case when $\Lambda$ is a block of a group ring).
\begin{defi}[Symmetrizing Form]
	Assume in this definition that the Wedderburn-decomposition of $A$ is given by
	\begin{equation}
		A \iso \bigoplus_{i=1}^n D_i^{d_i\times d_i}
	\end{equation}
	for certain skew-fields $D_i$  (finite-dimensional over $K$) and certain numbers $d_i$. 
	Let $\eps_i\in Z(A)$ (for $i\in \{1,\ldots,n\}$) be the central primitive idempotent belonging to the Wedderburn-component
	$D_i^{d_i\times d_i}$.
	For any element $u\in Z(A)^\times\iso \bigoplus_i Z(D_i)^\times$ define the non-degenerate associative symmetric bilinear form
	\begin{equation}
		T_u: \ A \times A \rightarrow K: \ (a,b) \mapsto \sum_{i=1}^n \operatorname{ tr.}_{Z(D_i)/K}\trred_{D_i^{d_i\times d_i} / Z(D_i)} (\eps_i \cdot u \cdot a \cdot b)
	\end{equation}
	Here $\operatorname{ tr.}_{Z(D_i)/K}: \ Z(D_i) \rightarrow K$ denotes the usual trace for 
	field-extensions.

	We call a full $\OO$-lattice $L\subset A$ \emph{self-dual} with respect to $T_u$ if it is equal to its \emph{dual lattice}
	$L^\sharp := \{  a \in A \ | \ T_u(a,L) \subseteq \OO\}$.
	If $\Lambda$ is self-dual with respect to $T_u$, then we call $T_u$ a \emph{symmetrizing form} for 
	$\Lambda$, and $u$ a \emph{symmetrizing element}. 
\end{defi}

\begin{remark}
	\begin{enumerate}
	 \item Any non-degenerate  symmetric and associative $K$-bilinear form on $A$ is equal to $T_u(-,=)$ for some $u\in Z(A)^\times$.
		This follows fairly easily from the structure theory of finite-dimensional semisimple algebras.
	 \item We sometimes write $T_u(a)$ (where $a\in A$) instead of $T_u(a,1)$.
	\end{enumerate}
\end{remark}

\begin{thm}[Transfer of the Symmetrizing Form]\label{thm_transfer_symm_form}
	Let $\Lambda$ be symmetric, and let $T\in \mathcal{K}^b(\projC_\Lambda)$ be a tilting complex.
	Set $\Gamma = \End_{\mathcal{D}^b(\Lambda)}(T)$, and $B=\End_{\mathcal{D}^b(A)}(K\otimes T)$.
	Identify
	\begin{equation}
		Z(A) = \bigoplus_{j=1}^n Z(\End_A(V_j)) = Z(B)
	\end{equation}
	Let $u = (u_1,\ldots,u_n) \in Z(A)^\times$ such that $\Lambda$ is self-dual with respect to the
	trace bilinear form $T_u: A \times A \rightarrow K$ induced by $u$. Then $\Gamma$ is self-dual
	with respect to the trace bilinear form $T_{\tilde{u}}: B\times B \rightarrow K$, where
	\begin{equation}
		\tilde{u} = ((-1)^{\eps(1)}\cdot u_1, \ldots, (-1)^{\eps(n)}\cdot u_n) \in Z(B)^\times
	\end{equation}
	where $\eps$ is as defined in Definition \ref{defi_data_tilting}.
\begin{proof}
	Let $\hat{u} = (\hat{u}_1,\ldots,\hat{u}_n)\in Z(B)$ be an element such that $\Gamma$ is actually self-dual with respect to $T_{\hat{u}}$. Then the $p$-valuations of the $\hat{u}_i$ are in fact independent of the particular choice of $\hat{u}$,
	since the coset $\hat{u} \cdot Z(\Gamma)^\times \in Z(B)^\times/Z(\Gamma)^\times$ is. Furthermore, the $u'\in Z(B)^\times$ such that $\Gamma$ is integral with respect to
    $T_{u'}$ are precisely the elements of $\hat u \cdot (Z(\Gamma) \cap Z(B)^\times)$. An element of $\hat u \cdot Z(\Gamma)\cap Z(B)^\times$ lies in
    $\hat u \cdot Z(\Gamma)^\times$ if and only in $\nu_p (u'_i) = \nu_p(\hat u_i)$ for all $i$ (all of those assertions are elementary). 
	Now assume we had shown that $\Gamma$ is integral with respect to $T_{\tilde u}$. Then we have 
    $\tilde u \in \hat u \cdot (Z(\Gamma) \cap Z(B)^\times)$. Thus $\nu_p (\tilde u_i) \geq \nu_p(\hat u_i)$ for all $i$, and equality for 
    all $i$ holds if and only if $\tilde u \in \hat u \cdot Z(\Gamma)^\times$, that is, if $\Gamma$ is self-dual with respect to
	$T_{\tilde u}$. So we have seen (up to the assumption above that we have yet to prove) 
    that if $\Lambda$ is self-dual with respect $T_{u}$ and $\Gamma$ is self-dual with respect to 
	$T_{\hat u}$, then $\nu_p(u_i) \geq \nu_p(\hat u_i)$, and, by swapping the roles of $\Lambda$ and $\Gamma$, also
    $\nu_p(\hat u_i) \geq \nu_p(u_i)$. In conclusion, we have $\nu_p(\tilde u_i) = \nu_p(u_i) = \nu_p(\hat u_i)$ for all $i$, which,
    by the above considerations, implies that $\Gamma$ is self-dual with respect to $T_{\tilde u}$.

    So far we have reduced the problem to showing that
$\Gamma$ is integral with respect to 
	$T_{\tilde{u}}$, which we will do now.
	So let $\varphi \in \End_{\mathcal{D}^b(\Lambda)}(T)$ (and fix a representative in $\End_{\mathcal C^b(\projC_\Lambda)}(T)$). Then 
	$\varphi^i$ induces an endomorphism of $T^i$, and we can decompose the $A$-module $K\otimes T^i$ as follows
	\begin{equation}
	K\otimes T^i = \underbrace{H^i(K\otimes T)}_{=: H^i} \oplus 
	\underbrace{\Im(\id_K\otimes d^{i-1})}_{=: Z^i_-} \oplus \underbrace{K\otimes T^i / \Ker(\id_K\otimes d^{i})}_{=:Z^i_+}
	\end{equation} 
	Define $\pi_{H^i}$, $\pi_{Z^i_-}$ and $\pi_{Z^i_+}$ to be the corresponding projections.
	Define $B^i := \End_A(K\otimes T^i)$, $B^i_{H} := \pi_{H^i} B^i \pi_{H^i} = \End_A(H^i)$,
	$B^i_+ := \pi_{Z^i_+} B^i \pi_{Z^i_+} = \End_A(Z^i_+)$ and $B^i_- := \pi_{Z^i_-} B^i \pi_{Z^i_-} = \End_A(Z^i_-)$. Now we have
	\begin{equation}
		\begin{array}{rcl}\displaystyle
		&&\displaystyle\sum_{i} (-1)^i\cdot T_{1_{B^i}\cdot u} (\varphi^i)\\ &=& \displaystyle \sum_i T_{\pi_{H^i}\cdot \tilde{u}}(\pi_{H^i}\varphi^i\pi_{H^i})
		+ (-1)^i\cdot T_{\pi_{Z^i_+}\cdot u}(\pi_{Z^i_+}\varphi^i\pi_{Z^i_+}) + (-1)^i \cdot T_{\pi_{Z^i_-} \cdot u}(\pi_{Z^i_-}\varphi^i\pi_{Z^i_-}) \\ &\stackrel{(*)}{=}& 
		\displaystyle\sum_i T_{\pi_{H^i}\cdot \tilde{u}}(\pi_{H^i}\varphi^i\pi_{H^i}) \stackrel{(**)}{=} \displaystyle T_{\tilde{u}}(\varphi)
		\end{array}
	\end{equation}
	Here $(*)$ holds because 
	\begin{equation}
		T_{\pi_{Z^i_+}\cdot u}(\pi_{Z^i_+}\varphi^i\pi_{Z^i_+})=T_{\pi_{Z^{i+1}_-}\cdot u}(\pi_{Z^{i+1}_-}\varphi^{i+1}\pi_{Z^{i+1}_-})
	\end{equation} as $\varphi$ is a map of
	chain complexes. The equality $(**)$ holds in fact just by definition, as we have identified $\bigoplus \End_A(H^i) = B$. The left side is trivially integral, as $\varphi^i \in \End_\Lambda (T^i)$, and $\End_\Lambda (T^i)$
	is a self-dual (and so in particular integral) lattice in $B^i$ with respect to $T_{1_{B^i}\cdot u}$. Hence the right side is also integral. So $\Gamma$ is indeed integral with respect to $T_{\tilde{u}}$. This concludes the proof.
\end{proof}
\end{thm}

\section{Partitioning $\Lifts(\overline\Lambda)$ by Rational Conditions}

Now we continue with what we started in Section \ref{section_corr_lifts}.
We want to define ``rational conditions'' on lifts that behave well under the map $\Phi_X$, that is, conditions such that $\Phi_X$ restricts to a bijective map between the lifts of $\overline\Lambda$
that fulfill the given conditions and the lifts of $\overline\Gamma$ that fulfill certain corresponding conditions. Probably the simplest of those conditions is to demand that the $K$-span of $\Lambda$ shall be Morita-equivalent to a certain semisimple $K$-algebra $A$.
It follows from the previous section that $\Phi_X$ sends lifts of $\overline\Lambda$ with $K$-span Morita-equivalent to
$A$ to lifts of $\overline\Gamma$ with the same property (and $\Phi_X^{-1}$ does it the other way round).

Since it will make things easier for us, we first give a slightly non-standard definition of decomposition matrices (which is linked
to the usual definition via Brauer reciprocity, and coincides with the usual definition in the split case).
\begin{defi}[Decomposition matrix]\label{defi_decomp_mat}
	We define the decomposition matrix $D^\Lambda$ of an $\OO$-order $\Lambda$ with $K\otimes \Lambda$ semisimple to be the transposed of the
	matrix of the canonical map of Grothendieck groups
	\begin{equation}
		K_0(\projC_{k\otimes\Lambda}) \iso K_0(\projC_\Lambda) \rightarrow K_0(\modC_{K\otimes \Lambda})
	\end{equation}
	sending $[P]$ to $[K\otimes P]$ with respect to the bases consisting of projective indecomposable modules on the left and simple $K\otimes \Lambda$-modules on the right. 
	 We call this map the ``decomposition map''.
	 Note that the rows of $D^\Lambda$ may be thought of as being labeled by the
	central primitive idempotents in $Z(K\otimes \Lambda)$ (resp. Wedderburn-components of $K\otimes\Lambda$).
\end{defi}

\begin{thm}
	Let $\overline{\Lambda}$ and $\overline{\Gamma}$ be finite-dimensional $k$-algebras that are derived equivalent.
	Let the derived equivalence be afforded by the (one-sided) tilting complex $\overline T$, and let 
	$X$ be a two-sided tilting complex such that its inverse has restriction to $\mathcal D^b(\projC_{\overline\Lambda})$ isomorphic to 
	$\overline T$. Set $\Phi := \Pi \circ \Phi_X$. Define 
	\begin{equation}
		\Lifts_s(\overline\Lambda) := \{ (\Lambda,\varphi)\in\Lifts(\overline\Lambda) \ | \ K\otimes\Lambda \textrm{ is semisimple } \}	 
	\end{equation}
 	Then $\Phi_X$ induces a bijection
	\begin{equation}
		\Lifts_s(\overline\Lambda) \longleftrightarrow \Lifts_s(\overline\Gamma)
	\end{equation}
	The following holds:
	\begin{enumerate} 	
	\item[(i)] If $(\Lambda,\varphi), (\Lambda',\varphi')\in \Lifts(\overline{\Lambda})$ are two lifts with $Z(K\otimes \Lambda) \iso Z(K\otimes \Lambda')$, then 
	\begin{equation}
		Z(K\otimes \Phi(\Lambda,\varphi)) \iso Z(K\otimes \Phi(\Lambda',\varphi'))
	\end{equation}
	and every choice of an isomorphism $\gamma: Z(K\otimes \Lambda)\rightarrow Z(K\otimes \Lambda')$ gives rise to a (canonically defined) isomorphism $\Phi(\gamma): Z(K\otimes \Phi(\Lambda,\varphi))\rightarrow Z(K\otimes \Phi(\Lambda',\varphi'))$.
	\item[(ii)] If $(\Lambda,\varphi), (\Lambda',\varphi')\in \Lifts(\overline{\Lambda})$ are two lifts and  $\gamma: Z( \Lambda) \stackrel{\sim}{\rightarrow} Z(\Lambda')$ is an isomorphism, then $\Phi(\gamma): Z(\Phi(\Lambda,\varphi)) \rightarrow Z(\Phi(\Lambda',\varphi'))$ is well defined and an isomorphism as well.
	\item[(iii)] If $(\Lambda,\varphi), (\Lambda',\varphi')\in \Lifts_s(\overline{\Lambda})$ are two lifts, and
	$\gamma: Z(K\otimes \Lambda) \stackrel{\sim}{\rightarrow} Z(K\otimes \Lambda')$ is an isomorphism such that  $D^\Lambda = D^{\Lambda'}$
	up to permutation of columns
	(where rows are identified via $\gamma$), then 
	$D^{\Phi(\Lambda,\varphi)}=D^{\Phi(\Lambda',\varphi')}$ up to permutation of columns (where rows are identified via $\Phi(\gamma)$).
	\item[(iv)]  If $(\Lambda,\varphi), (\Lambda',\varphi')\in \Lifts_s(\overline{\Lambda})$ are two lifts with  $D^\Lambda = D^{\Lambda'}$
	up to permutation of rows and columns then
	$D^{\Phi(\Lambda,\varphi)} = D^{\Phi(\Lambda',\varphi')}$
	up to permutation of rows and columns.
	\end{enumerate}
\begin{proof}
	 The fact that $\Phi_X$ induces a bijection between $\Lifts_s(\overline \Lambda)$ and $\Lifts_s(\overline \Gamma)$ follows
	 from the last section.

	Let $(\Lambda,\varphi) \in \Lifts(\overline{\Lambda})$. Then $Z(K\otimes \Phi(\Lambda))$
	is naturally isomorphic to $K\otimes Z(\Phi(\Lambda))$. But there is an isomorphism 
	between $Z(\Lambda)$ and $Z(\Phi(\Lambda,\varphi))$ (letting $c\in Z(\Lambda)$ correspond to the endomorphism of the
	tilting complex that is given by  multiplication with $c$ in every degree). That proves (i), and shows how $\Phi(\gamma)$ should be defined. 
	The claim of (ii) also follows.

	To the proof of (iii): Let $T\in\mathcal C^b(\projC_\Lambda)$ be the
	 lift of $\overline T$ (we identify $k\otimes\Lambda$ and $\overline\Lambda$ via $\varphi$). Write $\overline{T} = \overline{T}_0 \oplus \overline{T}_1$
	such that $\mathcal{G}_{\overline{T}}(\overline{T}_0) \iso 0 \rightarrow \overline{P} \rightarrow 0$ for a projective indecomposable $\overline{\Gamma}$-module $\overline{P}$.
	By Remark \ref{remark_partial_tilting} there is a corresponding direct sum decomposition 
	$T = T_0 \oplus T_1$ and we will have $\mathcal{G}_{T}(T_0) \iso 0 \rightarrow P \rightarrow 0$, where $P$ is the unique projective indecomposable $\Phi(\Lambda,\varphi)$-module
	with $k\otimes P \iso \overline{P}$. Then take $e_P$ to be the endomorphism of $T$ inducing the identity 
	on $T_0$ and the zero map on $T_1$. Clearly this is a primitive idempotent in $\Phi(\Lambda,\varphi)$
	(which is just $\End_{\mathcal{D}^b(\Lambda)}(T)$, so this statement makes sense) with $e_P\Phi(\Lambda,\varphi) \iso P$. So the decomposition number associated to $P$ and the simple $K\otimes \Phi(\Lambda,\varphi)$-module corresponding to the simple $K\otimes \Lambda$-module $V_j$ (under the isomorphism of the centers) is just the $\End_{K\otimes \Lambda}(V_j)$-rank of the image of $e_P$ in $\End_{K\otimes \Lambda}(V_j)^{\delta_j\times\delta_j}$ under the map given in Lemma \ref{lemma_tilting_wedderburn}. On the other hand (due the way Lemma  \ref{lemma_tilting_wedderburn} was obtained) this is just the absolute value of the coefficient of $[V_j]\in K_0(\modC_{K\otimes\Lambda})$ in the image under the decomposition map of $\sum_i (-1)^{i} \cdot[T^i_0]\in K_0(\projC_{\Lambda})$. But due to the isomorphism 
	$K_0(\projC_\Lambda)\iso K_0(\projC_{\overline{\Lambda}})$ we can compute this coefficient, and hence the decomposition matrix of 
	$\Phi(\Lambda,\varphi)$, from the knowledge of a direct sum decomposition of $\overline{T}$ and the
	knowledge of the decomposition matrix of $\Lambda$ (since the latter determines the map $K_0(\projC_{\overline\Lambda})\longrightarrow K_0(\modC_{K\otimes\Lambda})$). Therefore, if the decomposition matrices of $\Lambda$ and $\Lambda'$ coincide, then so do the decomposition matrices of $\Phi(\Lambda,\varphi)$ and $\Phi(\Lambda',\varphi')$.
	This concludes the proof of (iii). The explicit formula for the decomposition matrix of $\Phi(\Lambda,\varphi)$
    we obtained above
	is in fact independent of the knowledge of $Z(K\otimes \Lambda)$. This implies (iv).
\end{proof}
\end{thm}

\begin{remark}
    The last theorem shows that the lifts $(\Lambda,\varphi)\in \Lifts(\overline\Lambda)$ that satisfy 
    certain conditions (as listed in the theorem) correspond via $\Phi_X$ to lifts $(\Gamma,\psi)\in\Lifts(\overline\Gamma)$
    that satisfy a corresponding set of conditions. We shall call these kinds of conditions on $\Lambda$ ``rational conditions''.
\end{remark}

\section{$2$-Blocks with Dihedral Defect Group}

In this section we specialize $K$ to be the $2$-adic completion of the maximal unramified extension of $\Q_2$ (so, in particular, $k$ will be algebraically closed). We fix a finite group $G$ and a block $\Lambda$ of $\OO G$ with dihedral defect group $D_{2^n}$ for some fixed $n\geq 3$ (we use the convention
where $|D_{2^n}|=2^n$). Set $A:=K\otimes\Lambda$ and $\overline{\Lambda} := k\otimes\Lambda$. For any $i\geq 2$ we denote by $\zeta_i$ a primitive $2^i$-th root of unity in $\bar{K}$ (that is, we fix a choice for each $i$). In what follows, by a 
``character'' we always mean an absolutely irreducible ordinary character with values in $\bar{K}$.

\subsection{Generalities}

\begin{lemma}[Facts from Number Theory]\label{lemma_facts_from_num_theo}
	\begin{enumerate}
	 \item[(i)] Define $K_i := K(\zeta_{i}+\zeta_{i}^{-1})$. $K_i/K$ is a field extension of degree $2^{i-2}$. Its Galois group
		is cyclic, and we denote by $\gamma_i$ one of its generators. 
		Hence the subfield lattice of $K_{i}$ is just a chain, and in fact equal to
		\begin{equation}
			K=K_2 \subset K_3 \subset \ldots \subset K_{i}
		\end{equation}
		We denote by $\OO_i$ the integral closure of
		$\OO$ in $K_i$. 
	 \item[(ii)] The field extension $K_i/K$ is totally ramified and the $2$-valuation of its discriminant
		is equal to $(i-1)\cdot 2^{i-2}-1$. 
	 \item[(iii)] If $G$ is any finite group, then $KG$ is isomorphic to a direct sum of matrix rings over fields (i. e., no non-commutative division algebras occur in the Wedderburn decomposition of $KG$). 
	\end{enumerate}
\begin{proof}\begin{enumerate}
              \item[(i)] This is elementary Galois theory. 
	      \item[(ii)] This is \cite[Theorem 1]{LiangDiscrim}. The result from that paper carries over to our situation without change, as the
    $2$-valuation of the discriminant of $K(\zeta_i+\zeta_i^{-1})/K$ equals the $2$-valuation
    of the discriminant of $\Q(\zeta_i+\zeta_i^{-1})/\Q$ due to both extensions having the same degree.
	      \item[(iii)] To see this let 
	 $D$ be a skew-field that occurs in the Wedderburn decomposition of $\Q_2G$, and denote 
	 the center of $D$ by $E$. Then by \cite[Corollary 31.10]{Reiner} the unique unramified extension $E'$
	 of $E$ of degree equal to the index of $D$ will split $D$. 
	 We may write $E'=E\cdot F$ for some unramified extension $F$ of $\Q_2$. Then $F\otimes_{\Q_2} D$ will be isomorphic to
	 a direct sum of matrix rings over $E'$. Since $F$ is contained in $K$, this proves the assertion.
	    \end{enumerate}
\end{proof}
\end{lemma}

\begin{thm}[Brauer]\label{thm_brauer_char}
\begin{enumerate}
 \item[(i)] There are precisely $2^{n-2}+3$  characters in $\Lambda$. Four of these characters
	have height zero, the rest has height one. See \cite[Theorem 1]{BrauerDihedral}.
 \item[(ii)] All characters of $\Lambda$ take values in $K_{n-1}$ (see \cite[Proposition (5A)]{BrauerDihedral}). There are exactly $5$ characters in $\Lambda$ with values in $K$. The remaining characters
	lie in families $F_r$ for $r=1,\ldots,n-3$, where each $F_r$ is a single $\Gal(K_{n-1}/K)$-conjugacy class
	of characters. Each $F_r$ consists of $2^r$ elements (see \cite[Theorem 3]{BrauerDihedral}). Together with Lemma \ref{lemma_facts_from_num_theo} (i) and elementary Galois theory the latter implies that a character in $F_r$ takes
	values in $K_{r+2}$.
\item[(iii)] The four characters of height zero in $\Lambda$ take values in $K$. See \cite[Theorem 4]{BrauerDihedral}.
\end{enumerate}
\end{thm}

Note that we may as well denote the one-element set containing the unique $K$-rational character of height one
by $F_0$, and use indices $r=0,\ldots,n-3$. The grouping of the characters into four height zero characters
and $n-2$ families $F_r$ of height one characters seems more natural in what follows.

\begin{corollary}
	From the above it follows immediately that $\Lambda$ is an $\OO$-order in
	\begin{equation}\label{eqn_dkjhjkd676dghdg}
		A = \bigoplus_{i=1}^4 K^{\delta_i \times \delta_i} \oplus \bigoplus_{r=0}^{n-3} K_{r+2}^{\delta'_r\times \delta'_r} \quad \textrm{ for certain }\delta_i,\delta_i' \in \Z_{>0}
	\end{equation}
	that is self-dual with respect to $T_u$, where $u=(u_1,u_2,u_3,u_4,\ldots,u_{n+2})\in Z(A)$ with
	$\nu_2(u_i)=-n$ for $i=1,\ldots,4$ and $\nu_2(u_i)=-n+1$ for $i>4$. Of course the analogous statement will
	hold for a basic order of $\Lambda$.
\begin{proof}
      Lemma \ref{lemma_facts_from_num_theo} (iii) implies that $K\otimes \Lambda$ is a direct sum of matrix rings over fields (and not merely skew-fields). Therefore $K\otimes\Lambda$ is Morita-equivalent to its center. From ordinary representation theory we know that given a finite group $G$ we have
      \begin{equation}\label{eqn_dhgezuvhevdg3783bj}
      	Z(KG) \iso \bigoplus_\chi K(\chi)
      \end{equation}
      where $\chi$ runs over representatives for all Galois conjugacy classes of absolutely irreducible characters of $G$ with values in the algebraic closure of $K$. Theorem \ref{thm_brauer_char} says that in a block of defect $D_{2^n}$ there are $n+2$ Galois conjugacy classes of characters (four $K$-valued characters of height zero and one conjugacy class of $K_{r+2}$-valued characters for each 
      $r=0,\ldots,n-3$). This shows that $K\otimes\Lambda\iso A$ (with $A$ as given in (\ref{eqn_dkjhjkd676dghdg}) for some choice of numbers $\delta_i$ and $\delta_i'$).
      As for the choice of $u$, note that the group ring $\OO G$ of a finite group $G$ is self-dual with respect to $T_u$ where $u$ is defined as follows:
      \begin{equation}
      	 u = \sum_{\chi\in \Irr_{\bar K}(G)} \frac{\chi(1)}{|G|} \eps_\chi \in Z(\Q G) \subset Z(KG)
      \end{equation}
      Here,  $\eps_\chi$ denotes the central primitive idempotent in  $Z(\bar K G)$ associated to $\chi$. The entry of this element $u$ in the 
      Wedderburn component of the right hand side of (\ref{eqn_dhgezuvhevdg3783bj}) associated to the absolutely irreducible character $\chi$ is just $\chi(1)/|G|$. Now the assertions on the heights of the characters in a dihedral block  made in Theorem \ref{thm_brauer_char} imply our assertion on the $p$-valuations of the $u_i$. The fact that this symmetrizing element $u$ carries over to a basic order may be seen as a consequence of 
      Theorem \ref{thm_transfer_symm_form} (since a Morita equivalence is a special case of a derived equivalence).
\end{proof}
\end{corollary}

\begin{thm}[Erdmann]
	The basic algebra of $\overline{\Lambda}$ is isomorphic to one of the algebras of dihedral type in
	the list given in the appendix of \cite{TameClass}. \emph{(Technically, this follows from
	\cite[Lemma IX.2.2]{TameClass} together with the fact that $\overline{\Lambda}$ is known to be of
	tame representation type and thus has to occur in the list.)}
\end{thm}

\begin{thm}[Holm and Linckelmann]
	\begin{enumerate}
	 \item[(i)] In Erdmann's classification, the algebras $\mathcal{D}(2A)^{\kappa,c}$ and 
		$\mathcal{D}(2B)^{\kappa,c}$, for any combination $\kappa =2^{n-2}\geq 1$ and $c\in\{0,1\}$, are
		derived equivalent. In particular, for fixed $n$, there are at most two derived equivalence classes
		of $2$-blocks over $k$ with defect group $D_{2^n}$ and two simple modules. See \cite{HolmDerEq}.
	 \item[(ii)] There is precisely one derived equivalence class of $2$-blocks over $k$ with defect group $D_{2^n}$ and three simple modules. See \cite[Theorem 1]{LinckelmannDihedral}.
	\end{enumerate}
\end{thm}

\subsection{Blocks with Two Simple Modules}
\newcommand{\DtwoB}{\ensuremath{\mathcal{D}(2B)^{\kappa,c}}}

Assume in this subsection that $\Lambda$ has precisely two isomorphism classes of simple modules.
We first assume that $\overline{\Lambda}$ is Morita-equivalent to $\DtwoB$ for some $c\in\{0,1\}$
and $\kappa = 2^{n-2}$ (the latter is implied by $\kappa+3=\dim_k Z(\DtwoB)=\dim_K Z(A) = 2^{n-2}+3$). Now let $\Lambda_0$ be a basic algebra of $\Lambda$. From \cite{TameClass} we know that 
$k\otimes {\Lambda}_0 \iso kQ/I$, where
\begin{equation}
	Q=
\xygraph{
!{<0cm,0cm>;<1cm,0cm>:<0cm,1cm>::}
!{(0,0) }*+{\bullet_{0}}="a"
!{(2,0) }*+{\bullet_{1}}="b"
"b" :@/_/_{\gamma} "a"
"a" :@/_/_{\beta} "b"
"a" :@(ld,lu)^{\alpha} "a"
"b" :@(rd,ru)_{\eta} "b"
}
\end{equation}
and 
\begin{equation}
I = \left\langle \beta\eta,\ \eta\gamma,\ \gamma\beta,\ \alpha^2-c\cdot \alpha\beta\gamma,\  \alpha\beta\gamma - \beta\gamma\alpha,\ \gamma\alpha\beta-\eta^\kappa\right\rangle
\end{equation}
We may assume the following rational structure on $\Lambda_0$
\begin{equation}\label{rat_cond_DtwoB}
	\begin{array}{ccccc}	
		\mathbf{Z(A)} & \mathbf{u} & \mathbf{0} & \mathbf{1} \\ \cmidrule[1.2pt]{1-4}
		K & u_1 & 1 & 0 \\
		K & u_1 & 1 & 0 \\
		K & u_2 & 1 & 1 \\
		K & u_2 & 1 & 1 \\\cline{1-4}
		K_{r+2} & u_3 & 0 & 1 & \ [\textrm{ exactly once for each } r=0,\ldots,n-3\ ]
	\end{array}
\end{equation}
where $u_1,u_2\in K$ have $2$-valuation $-n$ and $u_3\in K$ has $2$-valuation $-n+1$.
\begin{remark}
    We say that a lift $\Gamma$ of $k\otimes\Lambda_0$ satisfies the rational conditions given above if all of the following conditions hold:
    \begin{enumerate}
    \item[(i)] $K\otimes\Gamma$ is Morita equivalent to $K\oplus K \oplus K \oplus K \oplus \bigoplus_{r=0}^{n-3} K_{r+2}$ (so in particular
    $K\otimes\Gamma$ will be semisimple).
    \item[(ii)] The decomposition matrix of $\Gamma$ is as in (\ref{rat_cond_DtwoB}), where the individual rows pertain
    to the summand of the center that is given on the left of the table. 
    \item[(iii)] There exists some 
    $u = (u_1, u_1, u_2, u_2, u_3,\ldots,u_3)\in K^{n+2} \subseteq K\oplus K \oplus K \oplus K \oplus \bigoplus_{r=0}^{n-3} K_{r+2}$ with $\nu_2(u_1)=\nu_2(u_2)=-n$ and $\nu_2(u_3)=-n+1$ such that $\Gamma$ is self-dual with respect to $T_u$.
	\end{enumerate}
    We should probably also explain what we mean when we say that two lifts $\Gamma$ and $\Gamma'$ of $k\otimes\Lambda_0$
    subject to the above rational conditions have \emph{equal} center. The point is that the rows of the decomposition
    matrix of $\Gamma$ are canonically in bijection with the Wedderburn components of $Z(K\otimes\Gamma)$ (or, equivalently, central
    primitive idempotents in $K\otimes\Gamma$). Naturally we demand that there should be an isomorphism $\gamma: Z(\Gamma)\stackrel{\sim}{\longrightarrow}Z(\Gamma')$ such that if $\eps\in Z(K\otimes\Gamma)$ is a central primitive idempotent, then the rows
    in the respective decomposition matrices pertaining to $\eps$ respectively $(\id_K\otimes\gamma)(\eps)$ are equal (up to some fixed permutation of the columns).
\end{remark}

\begin{lemma}\label{lemma_lift_commutative}
	Let
	\begin{equation}
		\Gamma \subseteq \OO \oplus \OO \oplus \bigoplus_{r=0}^{n-3}\OO_{r+2}
	\end{equation}
	be a local $\OO$-order such that $k\otimes\Gamma$  is generated by a single nilpotent element $\eta$ (so, in particular, $k\otimes\Gamma = k[\eta]$). 
	Furthermore assume that $\Gamma$ is symmetric with respect to $T_u$, where 
	$u=(u_1,u_2,u_3\ldots,u_n)\in K \oplus K \oplus \bigoplus_{r=0}^{n-3} K_{r+2}$
	with $\nu_2(u_1)=\nu_2(u_2)=-n$ and $\nu_2(u_i)=-n+1$ for all $i>2$.
	Then for some $x\in k^\times$
	there exists a preimage $\hat \eta$ of $x\cdot \eta$ in $\Gamma$ of the form
	\begin{equation}
		(0, 4, \pi_0,\ldots,\pi_{n-3})
	\end{equation}
	where the $\pi_r$ are prime elements in the ring $\OO_{r+2}$.
\begin{proof}
	If $\hat\eta=(a,b,d_0\ldots,d_{n-3})$ is a preimage of $\eta$, then $a\in (2)_\OO$, and hence
	$\hat\eta - a\cdot(1,\ldots,1)$ is a preimage of $\eta$ as well. So we may assume without loss
	that $a=0$. 
 	Hence some non-zero scalar multiple of $\eta$ will have a preimage in $\Gamma$ of the following shape:
	\begin{equation}
		\hat{\eta} = (0,2^l, \pi_0,\ldots,\pi_{n-3}) \quad\textrm{ with } \pi_r \in \Jac(\OO_{r+2})
	\end{equation}
	Note that we do not know yet that the $\pi_r$ are prime elements in $\OO_{r+2}$. All we can say at this point 
	is $\nu_2(\pi_r)\geq 2^{-r}$ (because we know the ramification indices of the extensions $K_{r+2}/K$ to be
	$2^r$).
	The fact that $\Gamma$ is self-dual with respect to $u$ implies that
	\begin{equation} \label{eqn_index_1}
		\nu_2 \left( \left[\OO_{n-1}^{2^{n-2}+1} : \OO_{n-1}\otimes_\OO \Gamma\right] \right) 
		= \frac{1}{2}\left({2n+(2^{n-2}-1)(n-1)}\right)
	\end{equation}
	Here, for two $\OO_{n-1}$-lattices $N \subseteq M$ such that $M/N$ is a torsion module, 
	 we denote by $[M:N]$ the product of all elementary divisors of $M/N$ (of course, this is only well-defied up to units).
	The left-hand side of the above equation is  equal to the $2$-valuation of the determinant of the 
	$(2^{n-2}+1)\times(2^{n-2}+1)$ Vandermonde matrix $M(S)$ associated to the values 
	\begin{equation}S = \{s_1,\ldots,s_{2^{n-2}+1}\} := \{ 0, 2^l, \pi_r^{\alpha_r} \ | \
	r=0,\ldots,n-3, \ \alpha_r \in \Gal(K_{r+2}/K) \}
	\end{equation}
	But the factorization (note that we fix an arbitrary total ordering on the Galois groups $\Gal(K_i/K)$)
	\begin{equation}
		\begin{array}{ll}
		 \displaystyle\prod_{i > j} (s_i - s_j) =&\displaystyle \pm 2^l\cdot\left(\prod_{r=0}^{n-3}\prod_{\alpha \in \Gal(K_{r+2}/K)} \pi_r^{\alpha}\right)
\cdot\left(\prod_{r=0}^{n-3}\prod_{\alpha \in \Gal(K_{r+2}/K)} (2^l-\pi_r^{\alpha})\right) \\ &
 \displaystyle\cdot \prod_{r=0}^{n-3} \left(
	\left(\prod_{q = r+1}^{n-3} \prod_{\alpha\in\Gal(K_{r+2}/K)} \prod_{\beta\in\Gal(K_{q+2}/K)}(\pi_r^\alpha -
	 \pi_q^\beta) \right)\right. \\&\displaystyle\cdot\left.
	\prod_{\alpha>\beta\in\Gal(K_{r+2}/K)} (\pi_r^\alpha - \pi_r^\beta) 
\right)
		\end{array}
	\end{equation}
	of $\det M(S)$ yields the following estimate of its $2$-valuation:
	\begin{equation}\label{eqn_index_2}
		\begin{array}{rl}
		&\nu_2(\det M(S))\\ \geq& \displaystyle l + \sum_{r=0}^{n-3}\frac{1}{2^r}\cdot 2^r + \sum_{r=0}^{n-3}\frac{1}{2^r}\cdot 2^r
		+ \sum_{r=0}^{n-3}\left( \sum_{q=r+1}^{n-3} (2^{r+q}\cdot\frac{1}{2^q}) + \frac{1}{2}\nu_2 \discrim_K (K_{r+2}) \right) \\
		 = & \displaystyle l+2(n-2)+\sum_{r=0}^{n-3}\left( (n-3-r)\cdot 2^r +\frac{1}{2}\left(
		(r+1)\cdot 2^r-1
		\right) \right) \\
		=& \displaystyle \frac{1}{2} n + \frac{1}{8} 2^n n + l - \frac{1}{8}2^n -\frac{3}{2}
		\end{array}
	\end{equation}
	Here we used that for any $x \in \Jac(\mathcal \OO_i)$ we have 
    \begin{equation}
	   \begin{array}{rcll}\displaystyle
	   \nu_2 \prod_{\alpha>\beta \in \Gal(K_i/K)} (x^\alpha - x^\beta) &=& \nu_2\left(\left[\OO_{n-1}^{2^{i-2}} : \OO_{n-1}\otimes_\OO \OO[x]\right]\right) \\
	   &\geq& \displaystyle\nu_2\left(\left[\OO_{n-1}^{2^{i-2}} : \OO_{n-1}\otimes_\OO \OO_i\right]\right)  = \frac{1}{2}\nu_2(\discrim_K (K_i))
	   \end{array}
	 \end{equation}

	Now the right hand side of (\ref{eqn_index_1}) has to be greater than or equal to the right hand side of (\ref{eqn_index_2}). This implies
	$l \leq 2$. On the other hand, the assumptions on $u$ would imply that $\nu_2( T_u(\hat{\eta})) < 0$ if $l\leq 1$, which is of course impossible.
	Hence $l=2$, and in particular the ``$\geq$'' in (\ref{eqn_index_2}) is really an equality, which is easily seen to be equivalent to
	$\nu_2 (\pi_r) = 2^{-r}$ for all $r=0,\ldots,n-3$.
\end{proof}
\end{lemma}

\begin{thm}\label{thm_unique_lift_dtwob}
	If $\Gamma,\Gamma' \in \mathfrak{L}(\DtwoB)$ (where $\kappa=2^{n-2}$) satisfy the rational conditions 
	stated in (\ref{rat_cond_DtwoB}) and $Z(\Gamma)= Z(\Gamma')$, then $\Gamma \iso \Gamma'$.
	Furthermore, the existence of such a lift implies $c=0$.
\begin{proof}
	Our general approach is to determine the structure of $\Gamma$ up to some parameters, and then 
	conclude that these parameters are determined by the knowledge of $Z(\Gamma)$.
	We assume (without loss) that
	\begin{equation}\label{eqn_iuehiudeh3egzge3hjhe}
		\Gamma \subseteq \OO \oplus \OO \oplus \OO^{2\times 2} \oplus \OO^{2\times 2} \oplus \bigoplus_{r=0}^{n-3} \OO_{r+2}
	\end{equation}
	Choose lifts $\hat{e}_0$ and $\hat{e}_1$ in $\Gamma$ of the idempotents $e_0$ and $e_1$ in $\DtwoB$. Assume without 
	loss that these idempotents $\hat e_0$ and $\hat e_1$ are diagonal in each direct summand on the right hand side of (\ref{eqn_iuehiudeh3egzge3hjhe}) (this is of course only a non-trivial condition in the two summands which $2\times 2$-matrix rings), and identify in the obvious way
	\begin{equation}
		\begin{array}{cc}
		\displaystyle\Gamma_{00} := \hat{e}_0 \Gamma \hat{e}_0 \subseteq \OO\oplus\OO\oplus\OO\oplus\OO &
		\displaystyle\Gamma_{11} := \hat{e}_1 \Gamma \hat{e}_1 \subseteq \OO\oplus\OO\oplus \bigoplus_{r=0}^{n-3} \OO_{r+2} \\ \\
		\displaystyle\Gamma_{10} := \hat{e}_1 \Gamma \hat{e}_0 \subseteq \OO\oplus\OO &
		\displaystyle\Gamma_{01} := \hat{e}_0 \Gamma \hat{e}_1 \subseteq \OO\oplus\OO 
		\end{array}
	\end{equation}
	We first look at $\Gamma_{11}$. Note that $e_1\DtwoB e_1 \iso k[\eta]$, and therefore
	Lemma \ref{lemma_lift_commutative} tells us that there is a lift $\hat{\eta}\in \Gamma_{11}$ of
	some non-zero scalar multiple of $\eta$
	of the form $(0,4,\pi_0,\ldots,\pi_{n-3})$.

	Now we consider $\Gamma_{00}$. We may assume without loss that $\Gamma_{00}$ is equal to the row space of
	\begin{equation}
		\left[
			\begin{array}{cccc}
				1&1&1&1 \\
				0 & 2^a & x & y \\
				0 & 0 & 2^b & z \\
				0 & 0 & 0 & 2^n
			\end{array}
		\right]
		\quad \textrm{ for certain } a,b\in\Z_{>0}\textrm{ and }x,y,z \in (2)_\OO
	\end{equation}
	We may furthermore assume without loss that $\hat{\alpha}:=[0,2^a,x,y]$ is 
	a lift of a (non-zero) scalar multiple of $\alpha$. To see this first note that
	$\hat \alpha \notin \Gamma_{01}\cdot \Gamma_{10}+2\cdot \Gamma_{00}$, and therefore the image
	of $\hat\alpha$ in $\DtwoB$ will be of the form $c_1\cdot \alpha+c_2\cdot\beta\gamma+c_3\cdot\alpha\beta\gamma$
	with $c_1,c_2,c_3\in k$ and $c_1\neq 0$. For all $c_1,c_2\in k$ there is an automorphism of $\DtwoB$ with
	$\alpha\mapsto \alpha+c_1\cdot \beta\gamma+c_2\cdot \alpha\beta\gamma$, $\beta\mapsto\beta$, $\gamma\mapsto\gamma$
	and $\eta\mapsto\eta$ (to verify this just plug the right hand sides into the defining relations of \DtwoB).
	Thus we may replace $\alpha$ by an appropriate multiple of the image of $\hat\alpha$ in $\DtwoB$.

	Next we look at the trace form $T_u$ to get some restrictions on the parameters (by ``$\sim$'' we mean ``equal up to units in $\OO$''):
	\begin{equation}
		\begin{array}{ccl}
			T_u([1,1,1,1]) \sim 2^{-n}\cdot (2 +2\cdot \frac{u_2}{u_1}) \stackrel{!}{\in} \OO & \rimply & \frac{u_1}{u_2} \equiv -1 \mod (2^{n-1}) \\
			T_u([0,0,2^b,z]) \sim 2^{-n}\cdot (2^b+z) \stackrel{!}{\in} \OO & \rimply & z \equiv -2^b \mod (2^n) \\
			& \stackrel{\textrm{w.l.o.g.}}{\rimply}& z=-2^b \\
			T_u(\hat{\alpha}) = 2^{-n}\cdot (2^a +(x+y)\cdot \frac{u_2}{u_1}) \stackrel{!}{\in} \OO
			& \rimply&  x+y \equiv -\frac{u_1}{u_2}\cdot 2^a \mod (2^n) \\
			&\stackrel{\textrm{w.l.o.g.}}{\rimply}& x = 2^a - y
		\end{array}
	\end{equation}
	
	Now let $\hat{\gamma} \in \Gamma_{10}$ and $\hat{\beta}\in\Gamma_{01}$ be lifts of non-zero scalar multiples of $\gamma$ and $\beta$ such that
	\begin{equation}
		\hat{\beta}\cdot\hat{\gamma} = [0,0,2^b,-2^b] + \xi \cdot [0,0,0,2^n] \quad \textrm{ for some $\xi \in \OO$}
	\end{equation}
	Then we have
	\begin{equation}
		\hat{\gamma}\cdot\hat{\beta} = [2^b,-2^b+\xi\cdot 2^n,0,\ldots,0] \in \Gamma_{11}
	\end{equation}
	Since $\beta\eta=0$ we have $\frac{1}{2} \cdot \hat{\gamma}\cdot \hat{\beta}\cdot \hat{\eta} \in \Gamma$, and thus
	\begin{equation}
		T_u\left(\frac{1}{2} \cdot \hat{\gamma}\cdot \hat{\beta}\cdot \hat{\eta}\right)
		= u_2\cdot (-2^{b+1} + \xi 2^{n+1}) \sim 2^{b-n+1} \stackrel{!}{\in} \OO \quad \rimply b \geq n-1 
	\end{equation}
	But $a+b=n$ and $a,b$ are both strictly greater than zero. This implies $b=n-1$ and $a=1$.
	To summarize: At this point 
	we know that $\Gamma_{00}$ is equal to the row space of
	\begin{equation}
		\left[
			\begin{array}{cccc}
				1&1&1&1 \\
				0 & 2 & x & 2-x \\
				0 & 0 & 2^{n-1} & 2^{n-1} \\
				0 & 0 & 0 & 2^n
			\end{array}
		\right]
		\quad \textrm{ for some $x \in (2)_\OO$ }
	\end{equation}
	Note that (for large $n$) this row space will not be multiplicatively closed for all values of $x$. So
	this gives us a condition on $x$:
	\begin{equation}
		\hat{\alpha}^2-2\hat{\alpha}= [0,0,x^2-2x,x^2-2x] \stackrel{!}{\in} \left\langle [0,0,2^{n-1},2^{n-1}], [0,0,0,2^n]\right\rangle_\OO	 
	\end{equation}
	This is equivalent to $x^2 \equiv 2x \mod (2^{n-1})$, which in turn is equivalent to
	\begin{equation}
		x\equiv 0 \mod (2^{n-2}) \quad\textrm{or}\quad x\equiv 2 \mod (2^{n-2}) 
	\end{equation}
	For now let us assume $x\equiv 0 \mod (2^{n-2})$. Then $x=2^{n-2}\cdot \xi$ for some $\xi \in \OO$.
	But then
	\begin{equation}
		\hat{\alpha}^2 - 2\hat{\alpha} = \xi(2^{n-3}\xi-1)\cdot [0,0,2^{n-1},2^{n-1}]
	\end{equation}
	and
	\begin{equation}
		\hat{\alpha} \cdot \Gamma_{01}\cdot \Gamma_{10} + 2\cdot \Gamma_{00} \subseteq \left\langle[0,0,0,2^n]\right\rangle_{\OO} + 2\cdot \Gamma_{00} 
	\end{equation}
	Hence $\alpha^2$ and $\alpha\beta\gamma$ would be linearly independent over $k$ if $\xi(2^{n-3}\xi-1)\in\OO^\times$.
	The relation $\alpha^2-c\cdot \alpha\beta\gamma$ prohibits this though. 
	Therefore we must have $\xi(2^{n-3}\xi-1) \in (2)_\OO$, and thus $\alpha^2=0$. This implies the assertion that the existence of a lift implies $c=0$. 
	Furthermore, if $n > 3$,  the fact that  $\xi(2^{n-3}\xi-1) \in (2)_\OO$ implies $x\equiv 0 \mod (2^{n-1})$. 
	If $n=3$,  the fact that $\xi(2^{n-3}\xi-1) \in (2)_\OO$ implies that  either $x\equiv 0 \mod (2^{n-1})$ or 
	$x\equiv 2 \mod (2^{n-1})$. Had we started with the assumption $x\equiv 2 \mod (2^{n-2})$, we would in the same fashion  have 
	arrived at $x\equiv 2 \mod (2^{n-1})$ (again with the exception of $n=3$ where $x\equiv 0 \mod (2^{n-1})$ is
	also possible). Hence independent of our assumptions on $x$ it follows that either $x\equiv 0 \mod (2^{n-1})$
	or $x\equiv 2 \mod (2^{n-1})$, which means that $\Gamma_{00}$ is equal to the row space of
	\begin{equation}\label{eqn_djijhduhu33773dghhgd}
		\left[
			\begin{array}{cccc}
				1&1&1&1 \\
				0 & 2 & 0 & 2 \\
				0 & 0 & 2^{n-1} & 2^{n-1} \\
				0 & 0 & 0 & 2^n
			\end{array}
		\right]
		\quad\textrm{or}\quad
		\left[
			\begin{array}{cccc}
				1&1&1&1 \\
				0 & 2 & 2 & 0 \\
				0 & 0 & 2^{n-1} & 2^{n-1} \\
				0 & 0 & 0 & 2^n
			\end{array}
		\right]	      
	\end{equation}
	The row space of the second matrix is obtained from the
	row space of the first matrix by swapping the first two columns. This swapping of columns is induced by an automorphism of 
	$K\otimes \Gamma$. Hence we may assume that we are in the case where  $\Gamma_{00}$ is equal to the row space of the leftmost matrix in (\ref{eqn_djijhduhu33773dghhgd}). Note that the aforementioned
	automorphism which swaps the first two Wedderburn components of $Z(K\otimes \Gamma)$ might not fix $Z(\Gamma)$. This will however not matter to us since we only use that
	the projection of $Z(\Gamma)$ to all but the first two Wedderburn components is equal to the  projection
	of $Z(\Gamma')$ to all but the first two Wedderburn components (instead of $Z(\Gamma)=Z(\Gamma')$; in particular, we could have made a slightly stronger
	assertion in the statement of the theorem).

	Now if we project $\Gamma_{00}$ onto its last two Wedderburn components we get an order 
	$\Gamma_{00}':=\left\langle [1,1], [0,2]\right\rangle_\OO$. Clearly $\Gamma_{01}$ and
	$\Gamma_{10}$ are both $\Gamma_{00}'$-lattices with the natural action. However, $\Gamma_{00}'$ has only two
	non-isomorphic lattices $L$ with $K\otimes L \iso K\otimes \Gamma_{00}'$, namely $L_1=\OO\oplus \OO$ and 
	$L_2=\Gamma_{00}'$. Both of them are self-dual lattices in $K\otimes \Gamma_{00}'$. Assume 
	$\Gamma_{01}=L_1$ (if we assume $\Gamma_{01}\iso L_1$, we may as well assume equality, by means of conjugation).
	By self-duality of $\Gamma$, we would then have $\Gamma_{10} = 2^{n}\cdot L_1$, and hence
	$\Gamma_{01}\Gamma_{10}\subset \Jac^2(\Gamma_{00})$. But $\beta\gamma$ certainly is not contained in
	$\Jac^2(e_0\DtwoB e_0)$. Hence we have a contradiction. This implies (without loss)
	$\Gamma_{01}=L_2$ and $\Gamma_{10}=[2^{n-1}, -2^{n-1}]\cdot L_2$.

	All that is left to verify is that the choice of the $\pi_i$ in $\Gamma_{11}$
	can be reconstructed from $Z(\Gamma)$. But from our knowledge of 
	$\Gamma_{00}$ and $\Gamma_{11}$ we know that the following element is in $Z(\Gamma)$:
	\begin{equation}
		[0,4,0,4,\pi_0,\ldots,\pi_{n-3}] \in Z(\Gamma) \subset K \oplus K \oplus K \oplus K \oplus \bigoplus_{r=0}^{n-3} K_{r+2}
	\end{equation}
	Hence the natural homomorphism $Z(\Gamma) \rightarrow \Gamma_{11}$ is surjective. This concludes the proof.
\end{proof}
\end{thm}

\newcommand{\DtwoA}{\ensuremath{\mathcal{D}(2A)^{\kappa,c}}}
Now assume that $\overline{\Lambda}$ is Morita-equivalent to $\DtwoA$.
Then we may assume the following rational structure of $\Lambda_0$
\begin{equation}\label{rat_cond_DtwoA}
	\begin{array}{ccccc}	
		\mathbf{Z(A)} & \mathbf{u} & \mathbf{0} & \mathbf{1} \\ \cmidrule[1.2pt]{1-4}
		K & u_1 & 1 & 0 \\
		K & u_1 & 1 & 0 \\
		K & u_2 & 1 & 1 \\
		K & u_2 & 1 & 1 \\\cline{1-4}
		K_{r+2} & u_3 & 2 & 1 & \ [\textrm{ exactly once for each } r=0,\ldots,n-3\ ]
	\end{array}
\end{equation}
where $u_1,u_2\in K$ have $2$-valuation $-n$ and $u_3\in K$ has $2$-valuation $-n+1$.
We also know from \cite{HolmDerEq} that there is a tilting complex $\overline{T}\in \mathcal{K}^b(\projC_{\DtwoA})$ with $\End_{\mathcal{D}^b(\DtwoA)}(\overline{T}) \iso \DtwoB$ looking as follows:
\begin{equation}\label{eqn_jkljuirrryy33}
	\overline{T} = \left[0 \rightarrow P_1 \oplus P_1 \rightarrow P_0 \rightarrow 0\right] \oplus \left[ 0 \rightarrow P_1 \rightarrow 0 \rightarrow 0\right]
\end{equation}
Let $X$ be a two-sided tilting complex the inverse of which restricts to $\overline T$.
Then clearly $\Phi_{X}$ maps a lift of $\DtwoA$ satisfying the rational conditions (\ref{rat_cond_DtwoA}) to a lift of $\DtwoB$ satisfying the rational conditions (\ref{rat_cond_DtwoB}). Hence we get 
the following Corollary directly:
\begin{corollary}\label{corr_ceq0_impossible}
    If there is a $\Gamma \in \mathfrak L(\DtwoA)$ subject to the rational conditions 
	stated in (\ref{rat_cond_DtwoA}), then $c=0$.
	In particular,
	if $B$ is a $2$-block of $kG$ with defect group $D_{2^n}$ (where $n\geq 3$), and $B$ has exactly two simple modules, then $B$ is Morita-equivalent to either $\mathcal{D}(2A)^{\kappa,0}$ or $\mathcal{D}(2B)^{\kappa,0}$ with $\kappa = 2^{n-2}$.
\end{corollary}

\begin{corollary}\label{corr_unique_lift_dwoa}
	If $\Gamma,\Gamma' \in \mathfrak{L}(\DtwoA)$ (where $\kappa=2^{n-2}$) satisfy the rational conditions 
	stated in (\ref{rat_cond_DtwoA}) and $Z(\Gamma)= Z(\Gamma')$, then $\Gamma \iso \Gamma'$.
\begin{proof}
	 By Corollary \ref{corr_lifts_two_term} $\Phi_X$ induces a bijection between $\mathfrak{L}(\DtwoA)$ and $\mathfrak{L}(\DtwoB)$.
	Note that
	$\Phi_X$ maps the
	lifts of $\DtwoA$ satisfying rational conditions as in (\ref{rat_cond_DtwoB})
	to lifts of $\DtwoB$ satisfying rational conditions as in (\ref{rat_cond_DtwoA}).
	Hence our assertion follows from Theorem \ref{thm_unique_lift_dtwob}.
\end{proof}
\end{corollary}

\subsection{Explicit Computation of the Lifts}
In this section we will compute the unique lift of $\DtwoA$ explicitly (depending, of course, on a prescribed center). We know already that we may assume $c=0$.
\renewcommand{\DtwoB}{\ensuremath{\mathcal{D}(2B)^{\kappa,0}}}
\renewcommand{\DtwoA}{\ensuremath{\mathcal{D}(2A)^{\kappa,0}}}
Define a complex of $\DtwoB$-modules
\begin{equation}\label{eqn_skls88}
    \overline T := \underbrace{0 \rightarrow P_1 \oplus P_1 \stackrel{\left[\begin{array}{c}\scriptstyle\gamma\\ \scriptstyle\gamma\alpha\end{array}\right]}{\longrightarrow} P_0 \rightarrow 0}_{=: \overline{T}_0}\ \oplus \
    \underbrace{0\rightarrow P_1 \rightarrow 0 \rightarrow 0}_{=: \overline T_1}
\end{equation}
Here, for the sake of simplicity, we identify the generators of $\DtwoB$ with homomorphisms between projective indecomposables satisfying the same
relations as the original generators (as opposed to the opposite relations). We can do this since the algebra $\DtwoB$ is isomorphic 
to its opposite algebra (it even carries an involution).  

\begin{remark}
	The algebra $\DtwoA$ has $\Ext$-quiver
	\begin{equation}\label{rels_dtwoa}
	 	Q'=
\xygraph{
!{<0cm,0cm>;<1cm,0cm>:<0cm,1cm>::}
!{(0,0) }*+{\bullet_{0}}="a"
!{(2,0) }*+{\bullet_{1}}="b"
"b" :@/_/_{\gamma'} "a"
"a" :@/_/_{\beta'} "b"
"a" :@(ld,lu)^{\alpha'} "a"
}
	\end{equation}
	with ideal of relations
	\begin{equation}
		I'=\left\langle \gamma'\beta', \ \alpha'^2, \ (\alpha'\beta'\gamma')^\kappa - (\beta'\gamma'\alpha')^\kappa \right\rangle_{kQ'}
	\end{equation}
	where $\kappa=2^{n-2}$. Its Cartan matrix is
	\begin{equation}
		\left[ \begin{array}{cc} 4\kappa & 2 \kappa \\ 2\kappa & \kappa+1 \end{array} \right]
	\end{equation}
\end{remark}

\begin{lemma}
  $\overline T$ as defined in (\ref{eqn_skls88}) is a tilting complex with endomorphism ring $\DtwoA$.
\begin{proof}
    First note that $\gamma$ and $\gamma\alpha$ form a $k$-basis of $\Hom(P_1,P_0)$. and $\beta, \alpha\beta$ form a $k$-basis of
    $\Hom(P_0,P_1)$. Now let $\varphi = c_1 \cdot \beta + c_2 \cdot \alpha\beta \in\Hom(P_0,P_1)$. Then 
    \begin{equation}
	   \left[\begin{array}{c} \gamma \\ \gamma\alpha\end{array} \right] \cdot \varphi = 0 \lrimply \left[ \begin{array}{c} c_2 \cdot \gamma\alpha\beta \\ c_1 \cdot \gamma\alpha\beta \end{array}\right]= 0 \lrimply \varphi=0
    \end{equation}
    This implies $\Hom(\overline T_0, \overline T[-1])=0$ (already in $\mathcal C^b(\projC_{\DtwoA})$). $\Hom(\overline T_1, \overline T[-1])=0$
    is clear since in any degree at least one of these complexes is the zero module. Now assume 
    $\varphi = c_1 \cdot \gamma + c_2\cdot \gamma\alpha\in \Hom(P_1,P_0)$. Then clearly
    \begin{equation}
	   \varphi = \left[\begin{array}{cc} c_1 & c_2\end{array}\right] \cdot \left[ \begin{array}{c} \gamma \\ \gamma\alpha \end{array} \right]
    \end{equation}
    which implies that every chain map from $\overline T$ to $\overline T[1]$ is homotopic to zero.
  Furthermore $\overline{T}$ generates $\mathcal D^b(\DtwoB)$, since $P_1[1]$ is a summand of $\overline{T}$,
  and the mapping cone of the projection map $\overline T _0 \rightarrow \overline T_1 \oplus \overline T_1$ is
  isomorphic to $P_0[0]$. So we have seen that $\overline T$ is a tilting complex.

  Now we claim that the endomorphisms
  \begin{equation}
	\xymatrix{
	 P_1 \oplus P_1 \ar[d]_{\left[\begin{array}{cc} 0 & 1 \\ 0 & 0\end{array}\right]} \ar[rr]^{\left[ \begin{array}{c} \gamma\\\gamma\alpha \end{array} \right]} & & P_0   \ar[d]^\alpha \\
	 P_1\oplus P_1 \ar[rr] && P_0
	}
  \end{equation}
  (which we denote by $\alpha'$) and 
  \begin{equation}
   \xymatrix{
	 P_1\oplus P_1 \ar[rr]\ar[d]_{\left[\begin{array}{c} 0\\1 \end{array}\right]} && P_0\ar[d]  && P_1\ar[rr]\ar[d]_{\left[\begin{array}{cc} \eta & 0 \end{array}\right]} && \ar[d] 0 \\
	 P_1 \ar[rr] && 0 && P_1\oplus P_1\ar[rr] && P_0
    }
  \end{equation}
  (which we denote by $\beta'$ and $\gamma'$) together with the idempotent endomorphisms coming from the decomposition $\overline T = \overline T_0 \oplus \overline T_1$ (which we denote by 
  $e_0'$ and $e_1'$) generate the endomorphism ring of $\overline T$.
  To prove this, we determine the dimension of the subalgebra of $\End(\overline T)$ they generate.
  It should be noted that one can deduce from the shape of $\overline T$ and the Cartan matrix of $\DtwoB$ that
  the Cartan matrix of $\End(\overline T)$ is equal to that of $\DtwoA$.
  First look at the endomorphism ring of $\overline{T}_0$ in the category $\mathcal{C}^b(\DtwoB)$
  (which we identify as a subring of $\End(P_1\oplus P_1)\oplus \End(P_0)$). Here 
  $\alpha'$ and $\beta'\cdot\gamma'$ generate the subalgebra
  \begin{equation}\label{eqn_krlokriou33}
	\left\langle \left( \left[\begin{array}{cc} 1&0\\0&1\end{array}\right], 1\right) \right\rangle_k \oplus 
	\left\langle \left( \left[\begin{array}{cc} 0&1\\0&0\end{array}\right], \alpha\right) \right\rangle_k 
	\oplus \left( \left[ \begin{array}{cc} \eta k[\eta] & \eta k[\eta] \\ \eta k[\eta] & \eta k[\eta] \end{array} \right], 0\right)
  \end{equation}
  which has dimension $2+4\cdot 2^{n-2}$.
  The zero-homotopic endomorphisms generate the subspace
  \begin{equation}\label{eqn_iuibdgtdf2674}
	\left\langle \left( \left[\begin{array}{cc} \gamma\alpha\beta&0\\0&0\end{array}\right], \alpha\beta\gamma\right), \left( \left[\begin{array}{cc} 0&\gamma\alpha\beta\\0&0\end{array}\right],
	0 \right), \left( \left[\begin{array}{cc} 0&0\\\gamma\alpha\beta&0\end{array}\right],
	\beta\gamma \right), 
	\left( \left[\begin{array}{cc} 0&0\\0&\gamma\alpha\beta\end{array}\right],
	\alpha\beta\gamma \right)
	\right\rangle_k
  \end{equation}
  which has two-dimensional intersection with the vector space in (\ref{eqn_krlokriou33}). Hence the subalgebra of the endomorphism
  ring (in $\mathcal D^b(\DtwoB)$) of $\overline T_0$ generated by $\alpha'$ and $\beta'\cdot\gamma'$ is $2^n$-dimensional.
  Since we know the dimension of $\End(\overline T_0)$ to be $2^n$, it follows that $\alpha'$ and $\beta'\cdot\gamma'$ generate $\End(\overline T_0)$.
  
  With much less effort one can see that (in the category $\mathcal C^b(\DtwoB)$) we have
  $\Hom(\overline{T}_0, \overline T_1) \iso k[\eta]\oplus k[\eta]$, and 
  $\beta'$ generates this space as an $\End(\overline T_0)$-module. Similarly
   $\Hom(\overline{T}_1, \overline T_0) \iso \eta k[\eta]\oplus \eta k[\eta]$
  and $\gamma'$ generates this space as an $\End(\overline T_0)$-module. Furthermore
  $\gamma'\cdot\alpha'\cdot \beta' = \eta$ generates $\End(\overline T_1)=\End(P_1)$ as a $k$-algebra. The above considerations imply that $e_0',e_1',\alpha',\beta'$ and $\gamma'$ generate the endomorphism ring  (in $\mathcal D^b(\DtwoB)$) of $\overline T$ as a $k$-algebra. 

 Now one can easily verify that $\alpha',\beta'$ and $\gamma'$ satisfy the relations given in (\ref{rels_dtwoa}), and this is all we have to check, since we know that the endomorphism ring of $\overline T$ has the same dimension as $\DtwoA$.
\end{proof}
\end{lemma}

\begin{thm}\label{thm_lifts_dihedral_two_simple}
	 Define $K$-algebras $A$ and $B$ as follows:
	 \begin{equation}
		  A := K\oplus K \oplus K^{2\times 2} \oplus K^{2\times 2} \oplus \bigoplus_{r=0}^{n-3} K_{r+2}
		  \quad  B := K\oplus K \oplus K^{2\times 2} \oplus K^{2\times 2} \oplus \bigoplus_{r=0}^{n-3} K_{r+2}^{3\times 3} 
	 \end{equation}
	 Define idempotents $\hat e_0, \hat e_1 \in A$:
	 \begin{equation}
		\hat e_0 := \left( 1, 1, \left[ \begin{array}{cc} 1&0\\0&0\end{array} \right], \left[ \begin{array}{cc} 1&0\\0&0\end{array} \right],
		0, \ldots, 0\right) \quad \hat e_1 := 1_A - \hat e_0
	 \end{equation}	 
	 and define idempotents $\hat e'_0, \hat e_1'\in B$:
	 \begin{equation}
		\hat e_0' := \left( 1, 1, \left[ \begin{array}{cc} 1&0\\0&0\end{array} \right],
		\left[ \begin{array}{cc} 1&0\\0&0\end{array} \right], \left[ \begin{array}{ccc} 1&0&0\\0&1&0\\0&0&0\end{array} \right],
		\ldots, \left[ \begin{array}{ccc} 1&0&0\\0&1&0\\0&0&0\end{array} \right]\right) \quad \hat e_1' := 1_B - \hat e_0'
	 \end{equation}

	 Any lift $\Lambda$ of $\DtwoB$ subject to the rational conditions in (\ref{rat_cond_DtwoB}) 
	 is isomorphic to the $\OO$-order in $A$ generated by the idempotents $\hat e_0$, $\hat e_1$ and
	\begin{equation}\label{eqn_iojioio1}
		\begin{array}{rcllllllllllllllllll}
			\hat e_0 A \hat e_0\ni\hat\alpha &=& (0, & 2, & 0, & 2) \\
			\hat e_1 A \hat e_1\ni\hat\eta &=& &&(0,&4,&\pi_0,&\ldots,&\pi_{n-3}) \\
			\hat e_0 A \hat e_1\ni\hat\beta &=& && (1,&1) \\
			\hat e_1 A \hat e_0\ni\hat\gamma &=& && (2^{n-1},&2^{n-1})
		\end{array}
	\end{equation}
    for certain prime elements $\pi_i \in K_{i+2}$.
    Any lift $\Gamma$ of $\DtwoA$ subject to the rational conditions in (\ref{rat_cond_DtwoA}) 
	 is isomorphic to the $\OO$-order in $B$ generated by the idempotents $\hat e_0'$, $\hat e_1'$ and
	\begin{equation}\label{eqn_iojioio2}
		\begin{array}{rclllllllllllllllll}
			\hat e_0' B \hat e_0' \ni \hat\alpha' &=& \bigg(0,& 2,& 2, & 0, & \left[ \begin{array}{cc} 0&1\\0&2\end{array}\right],
			& \ldots, & \left[ \begin{array}{cc} 0&1\\0&2\end{array}\right] \bigg) \\\\
			\hat e_0' B \hat e_1' \ni \hat\beta' &=& &&\bigg( 1, & 1, & \left[ \begin{array}{c} 0 \\ 1\end{array} \right], & \ldots, &\left[ \begin{array}{c} 0 \\ 1\end{array} \right] \bigg) \\\\
			\hat e_1' B \hat e_0' \ni \hat\gamma' &=& &&\big( -2, & -2, & \left[ \begin{array}{cc} \pi_0 & -2 \end{array} \right], & \ldots, &\left[ \begin{array}{cc} \pi_{n-3} & -2 \end{array} \right] \big)
		\end{array}
	\end{equation}
    for certain prime elements $\pi_i \in K_{i+2}$.
    In particular, any block with dihedral defect group $D_{2^n}$ and two simple modules is isomorphic to an order
    of one of the above shapes.

    Furthermore, if $X$ is a two-sided tilting complex the inverse of which restricts to $\overline T$, the lifts of (\ref{eqn_iojioio1}) and (\ref{eqn_iojioio2}) with equal $\pi_i$ correspond to each other under the bijection $\Phi_X$.	
\begin{proof}
    We have already seen in the proof of Theorem \ref{thm_unique_lift_dtwob} that $\Lambda$ has to be as in (\ref{eqn_iojioio1}). We did however not see (and in general it is not true) that $\hat \alpha, \hat \beta, \hat \gamma$ and $\hat \eta$ may be 
   assumed to be lifts of the elements $\alpha,\beta,\gamma$ and $\eta$. What we did see is that 
   $\hat\alpha$ and $\hat\eta$ may be assumed to reduce to scalar multiples of $\alpha$ and $\eta$.
   Since we will need it below we now show that we may in fact assume that $\hat\alpha, \hat\gamma$ and $\hat\eta$ reduce to $\alpha,\gamma$ and $\eta$. To see that one simply verifies that for all $c_1,c_2,c_3,c_4\in k$ with $c_1,c_2,c_4\neq 0$ the following
  \begin{equation}
  	\DtwoB \longrightarrow \DtwoB: \ \alpha\mapsto c_1\alpha, \ \beta \mapsto \frac{c_4^\kappa}{c_1c_2}\beta
	+\frac{c_3c_4^\kappa}{c_1c_2^2}\alpha\beta, \  \gamma \mapsto c_2\gamma+c_3\gamma\alpha, \ \eta \mapsto c_4\eta
  \end{equation}
   defines an automorphism of $\DtwoB$.

    Now we show that $\Gamma$ as given in (\ref{eqn_iojioio2}) equals $\Phi_X(\Lambda)$.
	We choose
	\begin{equation}
	  T := \underbrace{0 \rightarrow \hat P_1 \oplus \hat P_1 \stackrel{\left[\begin{array}{c}\scriptstyle\hat\gamma\\ \scriptstyle\hat\gamma\hat\alpha\end{array}\right]}{\longrightarrow} \hat P_0 \rightarrow 0}_{=: {T}_0}\ \oplus \
    \underbrace{0\rightarrow \hat P_1 \rightarrow 0 \rightarrow 0}_{=: T_1}
	\end{equation}
	as a lift of $\overline{T}$ (where the $\hat P_i$ are the projective indecomposable
	$\Lambda$-modules). Now 
	\begin{equation}
	\xymatrix{
	 \hat P_1 \oplus \hat P_1 \ar[d]_{\left[\begin{array}{cc} 0 & 1 \\ 0 & 2\end{array}\right]} \ar[rr] & & \hat P_0   \ar[d]^{\hat\alpha} \\
	 \hat P_1\oplus \hat P_1 \ar[rr] && \hat P_0
	}		
	\end{equation}
	is a lift of $\alpha'$ (which we denote by $\hat\alpha'$), and
	\begin{equation}
		\xymatrix{
	 	 \hat P_1\oplus \hat P_1 \ar[rr]\ar[d]_{\left[\begin{array}{c} 0\\1 \end{array}\right]} && \hat P_0\ar[d]  && \hat P_1\ar[rr]\ar[d]_{\left[\begin{array}{cc} \hat \eta & -2 \end{array}\right]} && \ar[d] 0 \\
	 \hat P_1 \ar[rr] && 0 && \hat P_1\oplus \hat P_1\ar[rr] && \hat P_0
		}
	\end{equation}
	are lifts of $\beta'$ and $\gamma'$ (which we denote by $\hat \beta'$ and $\hat\gamma'$). We now have
	to calculate the action of those endomorphisms on homology. For that identify 
	$K \otimes P_0 \iso K \oplus K \oplus K \oplus K$ and 
	$K\otimes P_1 \iso K \oplus K \oplus \bigoplus_{r=0}^{n-3} K_{r+2}$. Only for the third and fourth Wedderburn-component we need to do any actual work. Choose $\left[ \begin{array}{cc}0&1\end{array}\right]$
	as a basis for the projection  of the kernel of the differential to the third Wedderburn-component, and
	$\left[ \begin{array}{cc}-2&1\end{array}\right]$ as a basis of the projection to the fourth Wedderburn-component. Now, for instance,
	\begin{equation}
		\left[ \begin{array}{cc}0&1\end{array}\right] \cdot \left[ \begin{array}{cc} 0&1\\0&2 \end{array} \right] = 2 \cdot \left[ \begin{array}{cc}0&1\end{array}\right] \quad \textrm{and} \quad\left[ \begin{array}{cc}-2&1\end{array}\right] \cdot \left[ \begin{array}{cc} 0&1\\0&2 \end{array} \right] = 0 \cdot \left[ \begin{array}{cc}-2&1\end{array}\right]
	\end{equation}
	which leads to the corresponding entries of $\hat\alpha'$ in the third and fourth Wedderburn-component.
\end{proof}
\end{thm}

\section*{Acknowledgments}
This work was supported by the priority program DFG SPP 1388 of the German Research Foundation (DFG).
I would also like to thank G. Nebe for pointing out to me the problem of lifting tame blocks, and K. Erdmann for helpful discussions regarding the subject matter of this paper. 

\addcontentsline{toc}{section}{References}
\bibliographystyle{plain}
\bibliography{refs}

\begin{thebibliography}{10}

\bibitem{BleherUnivDef}
F.~M. Bleher, G.~Llosent, and J.~B. Schaefer.
\newblock Universal deformation rings and dihedral blocks with two simple
  modules.
\newblock {\em J. Algebra}, 345(1):49 -- 71, 2011.

\bibitem{BrauerDihedral}
R.~Brauer.
\newblock {On 2-blocks with dihedral defect groups}.
\newblock {\em Symposia Math.}, 13:{366--394}, 1974.

\bibitem{TameClass}
K.~Erdmann.
\newblock {\em {Blocks of Tame Representation Type and Related Algebras}}.
\newblock Number 1428 in Lecture Notes in Mathematics. Springer, 1990.

\bibitem{HolmDerEq}
T.~Holm.
\newblock Derived equivalent tame blocks.
\newblock {\em Journal of Algebra}, 194(1):178 -- 200, 1997.

\bibitem{ZimmermannGeometryChainComplexes}
B.~Huisgen-Zimmermann and M.~Saor{\'{\i}}n.
\newblock Geometry of chain complexes and outer automorphisms under derived
  equivalence.
\newblock {\em Trans. Amer. Math. Soc.}, 353(12):4757--4777 (electronic), 2001.

\bibitem{JensenXuZDegenerations}
B.~T. Jensen, X.~Su, and A.~Zimmermann.
\newblock Degenerations for derived categories.
\newblock {\em J. Pure Appl. Algebra}, 198(1-3):281--295, 2005.

\bibitem{DerEq}
S.~K\"onig and A.~Zimmermann.
\newblock {\em {Derived Equivalences for Group Rings}}.
\newblock Number 1685 in Lecture Notes in Mathematics. Springer, 1998.

\bibitem{LiangDiscrim}
J.~Liang.
\newblock On the integral basis of the maximal real subfield of a cyclotomic
  field.
\newblock {\em J. reine angew. Math.}, 286-287:223--226, 1976.

\bibitem{LinckelmannDihedral}
M.~Linckelmann.
\newblock A derived equivalence for blocks with dihedral defect groups.
\newblock {\em Journal of Algebra}, 164(1):244 -- 255, 1994.

\bibitem{Reiner}
I.~Reiner.
\newblock {\em {Maximal Orders}}.
\newblock Academic Press Inc., 1975.

\bibitem{RickardOne}
J.~Rickard.
\newblock Morita theory for derived categories.
\newblock {\em J. London Math. Soc.}, 39(2):{436--456}, 1989.

\bibitem{RickardDerEqDerFun}
J.~Rickard.
\newblock {Derived Equivalences as Derived Functors}.
\newblock {\em J. London Math. Soc.}, 43(2):37 -- 48, 1991.

\bibitem{RickardLiftTilting}
J.~Rickard.
\newblock Lifting theorems for tilting complexes.
\newblock {\em Journal of Algebra}, 142(2):383 -- 393, 1991.

\bibitem{RouquierAutomorphFrench}
R.~Rouquier.
\newblock { Groupes d'automorphismes et \'{e}quivalences stables ou
  d\'{e}riv\'{e}es}.
\newblock preprint.

\bibitem{RouquierICM2006}
R.~{Rouquier}.
\newblock {Derived equivalences and finite dimensional algebras}.
\newblock {\em Proc. ICM (Madrid 2006)}, 2:191--221, 2006.

\bibitem{ZimmTiltedSymmIsSymm}
A.~Zimmermann.
\newblock {Tilted symmetric orders are symmetric orders}.
\newblock {\em Archiv der Mathematik}, 73(1):15--17, 1999.

\end{thebibliography}
\end{document}